\documentclass[11pt]{article}
\usepackage{amsbsy,amsfonts,amssymb}
\usepackage{rotating}
\usepackage{amssymb}
\usepackage{amsmath}
\usepackage{amsfonts}
\usepackage{amsthm}

\newtheorem{theo}{Theorem}[section]
\newtheorem{lemma}[theo]{Lemma}
\newtheorem{prop}[theo]{Proposition}

\newtheorem{btheo}{Theorem}[theo]

\newtheorem{coll}[theo]{Corollary}

\theoremstyle{definition}
\newtheorem{example}[theo]{Example}

\newtheorem{remark}[theo]{Remark}

\newtheoremstyle{ntheo}{}{}{}{}{\bfseries}{}{6pt}{(\thmname{#1}\thmnumber{#2})\thmnote{ \textup{(#3)} }}
\theoremstyle{ntheo}
\newtheorem{ntheo}[theo]{}

\newcommand{\be}{\begin{equation}}
\newcommand{\ee}{\end{equation}}
\newcommand{\ba}{\begin{array}}
\newcommand{\ea}{\end{array}}
\newcommand{\m}{\frak m}
\newcommand{\n}{\frak n}

\begin{document}

\title{Commutative rings with toroidal zero-divisor graphs }

\author{ Hung-Jen Chiang-Hsieh\thanks{e-mail: hchiang@math.ccu.edu.tw (corresponding author)}
~and Hsin-Ju Wang\thanks{e-mail:  hjwang@math.ccu.edu.tw}\\
\small{Department of Mathematics,}\\
\small{ National Chung Cheng University,}\\
\small{ Chiayi 621, Taiwan}\\
\vspace{.1in}\\
Neal O. Smith\thanks{e-mail: nsmith12@aug.edu}\\
\small{Department of Mathematics,}\\
\small{Augusta State University,}\\
\small{Augusta GA 30904}\\
}

\date{}

\maketitle

\begin{abstract} \noindent Let
$R$ be a commutative ring and let $\Gamma (R)$ denote its
zero-divisor graph. We investigate the genus number of the compact
Riemann surface in which $\Gamma (R)$ can be embedded and explicitly
determine all finite commutative rings $R$ (up to isomorphism) such
that $\Gamma(R)$ is either toroidal or planar. 
\end{abstract}
\section*{Introduction} We assume that all rings are commutative with
identity. For a ring $R$, the zero-divisor graph of $R$, denoted
by $\Gamma(R)$, is the simple graph whose vertex set consists of
all nonzero zero-divisors of $R$. Two distinct vertices are joined
by an edge if and only if the product of the vertices is $0$.
Therefore $\Gamma(R)=\emptyset$ if and only if $R$ is an integral
domain. This definition was introduced by Anderson and Livingston
in \cite{al}. Recently, this subject has been extensively studied
in \cite{amy}, \cite{am}, \cite{afll}, \cite{als}, \cite{rb},
\cite{ds}, \cite{la}, \cite{li}, \cite{mu}, \cite{ns1},
\cite{ns2}, and \cite{wa}.

There are many known results concerning zero-divisor graphs.
Anderson and Livingston showed in \cite{al} that $\Gamma(R)$ is
always connected and $R$ is a finite ring or an integral domain if
and only if $\Gamma(R)$ is finite.  Mulay \cite{mu} showed that if
$\Gamma(R)$ contains a cycle, then $\Gamma(R)$ contains a 3-cycle
or a 4-cycle. Anderson, Frazier, Lauve and Livingston showed in
\cite{afll} that if $R$ and $S$ are finite reduced rings which are
not fields, then $R\simeq S$ if and only if $\Gamma(R)\simeq
\Gamma(S)$.

The main objective of topological graph theory is to embed a graph
into a surface. Simply stated, that is to draw a graph on a
surface so that no two edges cross one another.  The simplest case
of this problem is when the surface in question is the plane; if a
graph can be embedded in the plane, we say the graph is
\emph{planar}. There are many papers where planarity of
zero-divisor graphs has been discussed.  In \cite{afll}, Anderson
et al. the authors determined when $\mathbb{Z}_{n_1} \times \dots
\times \mathbb{Z}_{n_k}$ and $\mathbb{Z}_n[X]/(X^m)$ have planar
zero divisor graphs and posed the general question as to which
finite rings $R$ have $\Gamma(R)$ planar. It was shown in
\cite{amy} that if $R$ is a finite local ring such that
$\Gamma(R)$ has at least 33 vertices, then $\Gamma(R)$ is not
planar. In that paper, Akbari et al. conjectured that for any
local ring of cardinality 32 which is not a field, $\Gamma(R)$ is
not planar.

This conjecture was proved independently in two papers.  In
\cite{ns1}, Smith proved that if $R$ is a finite, commutative local
ring (not a field) with cardinality 28 or greater, then $\Gamma(R)$
is not planar. Also in that paper the author classifies precisely
those finite commutative rings $R$ for which $\Gamma(R)$ is planar.
Some of the methods used in that paper are similar in spirit to some
of the arguments in section three of this paper; thus, we recover
this listing of planar graphs as we work toward our main theorem. We
also refer the readers to \cite{rb}, in which Belshoff and Chapman
independently obtain some of the same results from \cite{ns1}, and
to \cite{re}, in which Redmond lists all zero-divisor graphs up to
14 vertices where some of them are among the same results as above.

Akbari's conjecture was also verified independently in \cite{wa}.
Moreover, in the same paper the author also found all finite rings
of the form $\mathbb{Z}_{p_1^{\alpha_1}}\times \cdots \times
\mathbb{Z}_{p_n^{\alpha_n}}$ or of the form
$\mathbb{Z}_n[x]/(x^m)$ whose zero-divisor graphs are planar or
can be embedded into a torus.  This motivates the work in this
paper.

To find all finite rings $R$ such that $\Gamma(R)$ has genus at
most one is the goal of this paper. Since a finite ring is
Artinian, it is a direct product of local Artinian rings. Thus, we
first consider the case of finite local rings. To motivate the
main theorems in section 3, we first discuss the genera of local
rings under some specific assumptions in section 2. Using the
Euler characteristic formula and a technique of deletion and
insertion, we are able to successfully exclude some cases of
higher genus.

In section 3, we consider case by case those local rings $(R,\m)$
with  $|R/\m|\leq 8$. From \cite[Theorem~3.6]{wa}, we know that
these are all the cases of interest.  By \cite[Lemma 3.1]{wa}, if
a finite ring $R$ has $|Spec(R)|\geq 5$ then
$\gamma(\Gamma(R))\geq 2$. According to this, it suffices to look
for the finite rings with at most 4 maximal ideals. We obtain a
complete list of all finite rings whose zero-divisor graphs have
genera at most one in this section and summarize them in four
tables at the very end of this paper.  Finally, we do a similar
analysis in the case where $R$ decomposes as a product of local
rings.

We should mention that as this paper was being submitted, the
authors became aware of a similar work. In \cite{wi}, Wickham
independently proved some of the same results that appear in
section three of this paper.

\section{Preliminaries}
In this section we briefly recap some notation, terminology, and
basic results from \cite{hf}, \cite{ry} and \cite{wa}.

A simple graph $G$ is an ordered pair of disjoint sets $(V, E)$,
such that $V=V(G)$ is the set of vertices of $G$ and $E=E(G)$ is
the set of edges of $G$. For $v\in V$, the {\it degree} of $v$,
denoted by $deg(v)$, is the number of edges of $G$ incident to
$v$. If $V'\subseteq V(G)$, we define $G-V'$ to be the subgraph of
$G$ obtained by deleting the vertices in $V'$ and all incident
edges.  Similarly, if $E'\subseteq E(G)$, then $G-E'$ is the
subgraph of $G$ obtained by deleting the edges in $E'$. For a
graph $G$, let $\widetilde{G}$ denote the subgraph $G-V'$ where
$V'=\{\,v\in V~|~deg(v)=1\}$.  We call this graph the {\it
reduction} of $G$.
\par
A graph in which each pair of distinct vertices is joined by an
edge is called a {\it complete graph}.  We use $K_n$ to denote a
complete graph with $n$ vertices. A {\it bipartite graph} $G$ is a
graph whose vertex set $V(G)$ can be partitioned into two subsets
$V_1$ and $V_2$.  The edge set of such a graph consists of
precisely those edges which join vertices in $V_1$ to vertices of
$V_2$. In particular, if $E(G)$ consists of all possible such
edges, then $G$ called a {\it complete bipartite graph} and
denoted by the symbol $K_{m,n}$ where $|V_1|=m$ and $|V_2|=n$.

\par
By a \emph{surface}, we mean a two dimensional real manifold, that
is a topological space such that each point has a neighborhood
homeomorphic to the open disc. It is well-known that every
orientable compact surface is homeomorphic to a sphere with $g$
handles. This number $g$ is called the genus of the surface. For
example, the genus of a sphere is 0 and the genus of a torus is 1.
A simple graph which can be drawn without crossings on a surface
of genus $g$ but not on one of genus $g-1$, is called a graph of
genus $g$. We say a {\it planar} graph is a graph of genus 0 and a
{\it toroidal} graph is a graph of genus 1. We use $\gamma(G)$ to
denote the genus of a graph $G$. The following two results (see
\cite[p.118]{hf}) about the genus of a complete graph and a
complete bipartite graph will be very useful in the subsequent
sections.
\begin{lemma}
\label{genuskn} Let $n \geq 3$.
 Then $\gamma(K_n)=\{\dfrac{1}{12}(n-3)(n-4)\}$, where
$\{x\}$ is the least integer that is greater than or equal to $x$.
In particular, $\gamma(K_n)=1$ if $n=5, 6,$ or $7$.
\end{lemma}
\begin{lemma}
\label{genusbi} $\gamma(K_{m,n})=\{\dfrac{1}{4}(m-2)(n-2)\}$,
where $\{x\}$ is the least integer that is greater than or equal
to $x$. In particular, $\gamma(K_{4,4})=\gamma(K_{3,n})=1$ if
$n=3, 4, 5,$ or $6$.
\end{lemma}

Suppose a connected graph $G$ is drawn on an orientable compact
surface $S_g$ and let $\#V_G, \#E_G,$ and $\#F_G$ denote the
number of vertices, edges, and faces of $G$ respectively.  The
well known Euler characteristic formula states that
$\#V_G-\#E_G+\#F_G=2-2g$, where $g$ is the genus of $S_g$.

We end this section with two remarks.
\begin{remark}
\label{reduced}   $\gamma(H)\leq\gamma(G)$ for all subgraphs $H$
of $G$; and $\gamma(\widetilde{G})=\gamma(G)$, where
$\widetilde{G}$ is the reduction of $G$.
\end{remark}

\begin{remark}
\label{rem1} The bipartite graph $K_{3,6}$ has $v=9$ vertices and
$e=18$ edges. By Lemma \,\ref{genusbi}, we see that
$\gamma(K_{3,6})=1$. Therefore, from the Euler characteristic
formula there are $f=9$ faces when drawing $K_{3,6}$ without
crossings on a torus. We note that the boundary of each face $F_i$
of $K_{3,6}$ is an even cycle with length $e_i\geq 4$. It follows
from the inequality $2e=\sum_{i=1}^{f}e_i\geq 4f$ that $e_i=4$ for
each $i$.  That is, all face boundaries are $4$-cycles. Moreover,
any two faces in $K_{3,6}$ have at most one boundary edge in
common.
\end{remark}

\section{Genera of some special rings}

In order to simplify the proof of our main result in the next
section, we now discuss several special rings with
$\gamma(\Gamma(R))\geq 2$ or $\gamma(\Gamma(R))=1$. In the sequel,
if $R$ is a ring then $Z(R)$ will denote the set of its
zero-divisors and we define $Z(R)^*=Z(R)-\{0\}$.
\begin{prop}
\label{examlem1} Let $(R, \m)$ be a local ring with $|R|=32$. If
$|R/\m|=2$ and $|\m^2|=4$, then $\gamma(\Gamma(R))\geq 2$.
\end{prop}
\begin{proof}
We have that $|\m|=16$ and $dim_{R/\m} \m/\m^2=2$. Since $\m^2\neq
\{0\}$, by Nakayama's lemma \cite[Prop. 2.6]{ama} we have
\mbox{$\m^3\subsetneq\m^2$}.  Thus $\m^3=\{0\}$ or $|\m^3|=2$. If
$\m^3=\{0\}$, then $|\m-\m^2|=12$ and $|\m^2-\m^3|=3$. This
implies that $K_{3,12}\subseteq \Gamma(R)$ and therefore
\mbox{$3=\gamma(K_{3,12})\leq\gamma(\Gamma(R))$} by
Lemma~\ref{genusbi} and Remark~\ref{reduced}. Hence, we may assume
$|\m^3|=2$, so that $\m^4=\{0\}$ and $dim_{R/\m}
\m^2/\m^3=dim_{R/\m} \m^3/\m^4=1$.

Since $dim_{R/\m} \m/\m^2=2$, it follows from \cite[Prop.2.8]{ama}
that $\m$ can be generated by two elements.  Write $\m=(x, y)$.
Assume first that $x^3=y^3=0$. If $x^2y\neq 0$, then $x^2, xy\notin
\m^3$, so that $\{x^2\}$ and $\{xy\}$ are both bases for
$\m^2/\m^3$.  It then follows that $xy-x^2\in \m^3$ as $|R/\m|=2$.
Therefore, we have $x^2y-x^3=0$ and $x^2y=0$, a contradiction. Thus,
$x^2y=0$ and $xy^2=0$. If $x^3=y^3=0$ we have $\m^3=\{0\}$, a
contradiction. Consequently, we conclude that either $x^3\neq 0$ or
$y^3\neq 0$. We may assume without loss that $x^3\neq 0$. Therefore,
$\{x^3\}$ is a basis for $\m^3/\m^4$ and $\{x^2\}$ is a basis for
$\m^2/\m^3$. If $xy\notin \m^3$, then $xy-x^2\in \m^3$ as
$dim_{R/\m} \m^2/\m^3=1$ and $|R/\m|=2$, so that $x(y-x)\in \m^3$.
Therefore, we may replace $y$ by $y-x$ and assume that $xy\in \m^3$.
Moreover, if $xy\neq 0$, then $xy-x^3=0$ as $dim_{R/\m} \m^3/\m^4=1$
and $|R/\m|=2$, so that $x(y-x^2)=0$. Hence we may replace $y$ by
$y-x^2$ and assume that $xy=0$. We observe that $y^2\in \m^3$. If
not, then $y^2-x^2\in \m^3$, so that $xy^2-x^3=0$.  It then follows
that $x^3=0$, a contradiction. Now, we have two cases to discuss, as
either $y^2=0$ or $y^2=x^3$.

\noindent \textbf{Case 1} : $y^2=0$. In this case, let $u_1=x^3$,
$u_2=y$, $u_3=y+x^3$, $v_1=x$, $v_2=x+x^2$, $v_3=x+x^3$,
$v_4=x+x^2+x^3$, $v_5=x^2$, $v_6=x^2+x^3$ and $v_7=x+y$.  Then
$u_i\cdot v_j=0$ for every $i, j$. Therefore $K_{3,7}\subseteq
\Gamma(R)$ and it follows that $\gamma(\Gamma(R))\geq
\gamma(K_{3,7})=2$.\\

\noindent \textbf{Case 2} : $y^2=x^3$. In this case, let
$u_1=x^3$, $u_2=y$, $u_3=y+x^3$, $v_1=x$, $v_2=x+x^2$,
$v_3=x+x^3$, $v_4=x+x^2+x^3$, $v_5=x^2$, $v_6=x^2+x^3$, $w_1=x+y$,
$w_2=x+y+x^2$, $w_3=x+y+x^3$, $w_4=x+y+x^2+x^3$, $w_5=y+x^2$ and
$w_6=y+x^2+x^3$. Observe that $u_i\cdot v_j=0$ for every $i, j$,
so that $K_{3,6}\subseteq \Gamma(R)$.  Therefore
$\gamma(\Gamma(R))\geq \gamma(K_{3,6})=1$. Write $G=\Gamma(R)$,
$G'=G-\{u_1u_2, u_1u_3, u_1w_5, u_1w_6\}$\footnote{If $u,v\in R$
such that $uv=0$ then we write $uv$ as an edge of $\Gamma(R)$.},
and $G''=G'-\{w_1, \dots, w_6\}$.  It is then easy to see that
$G''\simeq K_{3,6}$. Next, we proceed to prove $\gamma(G)\geq 2$
by a deletion and insertion argument.

Suppose that $\gamma(G)=1$. Since $1=\gamma(G'')\leq
\gamma(G')\leq \gamma(G)$, we get $\gamma(G')=1$. Since
$$V(G')=\{u_1, u_2, u_3, v_1, \dots, v_6, w_1, \dots, w_6\}$$ and
$$\ba{rl} E(G')=&\,\,\,\,\,\, \{\,u_iv_j\,|\,1\leq i\leq 3, 1\leq j\leq 6\,\}\,\cup\, \{\,u_1w_i\,|\,1\leq i\leq 4\}\\
                &\cup\,\, \{\,w_iw_j\,|\,1\leq i\leq 4, j=5, 6\}\,\cup\, \{\,v_5w_5, v_5w_6, v_6w_5,
                v_6w_6\,\},\ea $$
by the Euler characteristic formula there are 19 faces when
drawing $G'$ on a torus. Fix a representation of $G'$ and let
$\{F'_1, \dots, F'_{19}\}$ be the set of faces of $G'$
corresponding to this representation. We note that $G''\simeq
K_{3,6}$ and therefore this graph has 9 faces whose boundaries are
all 4-cycles (see Remark~\ref{rem1}). Write $F''_1, \dots, F''_9$
for the faces of $G''$ obtained by deleting $w_1, \dots, w_6$ and
all edges incident with $w_1, \dots, w_6$ from the
representation of $G'$. 
Then $\{F'_1, \dots, F'_{19}\}$ can be recovered by inserting
$w_1, \dots, w_6$ and all edges incident with $w_1, \dots, w_6$
into the representation corresponding to $\{F''_1, \dots,
F''_9\}$. Let $F''_{t_i}$ denote the face of $G''$ into which
$w_i$ is
inserted during the recovering process from $G''$ to $G'$. %
We note that $w_iw_j\in E(G')$ for $i=1,\ldots,4$ and $j=5,6$.
Therefore every $w_i$ should be inserted into the same face, say
$F''_m$, of $G''$ to avoid any crossings, i.e.,
$t_1=t_2=\cdots=t_6=m$. Moreover, since $u_1w_i\in E(G')$ for
$i=1,\dots,4$, $u_1$ is a vertex of the face $F''_m$. Write the
edges $e_i=u_1w_i$, $e_{i+4}=w_iw_5$ and $e_{i+8}=w_iw_6$ for
$i=1,\dots, 4$.  After inserting $w_1, \dots, w_5$ and $e_1,
\dots, e_8$ into $F''_m$ we obtain Figure 1 as below. However, it
is easy to see from Figure 1 that we can not insert $w_6$ and
$e_9, \dots, e_{12}$ into $F''_m$ without crossings, a
contradiction. Therefore, we may conclude that $\gamma(G)\geq 2$.
\end{proof}
\begin{picture}(200,130)(-60,-120)
\put(150,0){\line(1, -1){60}}\put(210,-60){\line(-1, -1){60}}
\put(150,-120){\line(-1,1){60}}\put(90,-60){\line(1,1){60}}

\put(147.52,-3.6){$\bullet$}\put(156,-3){$u_1$}%
\put(88,-63){$\bullet$}%
\put(207,-63){$\bullet$}%
\put(147.52,-122){$\bullet$}%
\put(147.52,-92.36){$\bullet$}\put(146,-102){$w_5$}\put(161,-60){$\bullet$}\put(136,-61){$\bullet$}
\put(114,-60){$\bullet$}\put(181,-60){$\bullet$}\put(165,-60){$w_3$}\put(126,-60){$w_2$}\put(103,-60){$w_1$}
\put(185,-60){$w_4$}

\put(150,-1){\line(1,-4){14}}\put(150,-1){\line(-1,-5){11}}
\put(150,-1){\line(-3,-5){33}}\put(150,-1){\line(3,-5){34}}

\put(150,-90){\line(2,5){12.5}}\put(150,-90){\line(1,1){34}}\put(150,-90){\line(-1,1){34}}
\put(150,-90){\line(-2,5){12.5}}
\put(125,-80){$e_5$}\put(165,-80){$e_8$}\put(125,-30){$e_1$}\put(165,-30){$e_4$}\put(138,-37){$e_2$}
\put(154,-37){$e_3$}\put(137,-72){$e_6$}\put(154,-72){$e_7$}

\put(110,-6){$F''_m:$} \put(130,-130){Figure~1}
\end{picture}

\begin{prop}
\label{examlem2} Let $(R, \m)$ be a local ring with $|R|=32$. If
$|R/\m|=2$ and $|\m^2|=2$, then $\gamma(\Gamma(R))\geq 2$.
\end{prop}
\begin{proof}
Note that $|\m|=16$, $\m^3=\{0\}$, $dim_{R/\m} \m/\m^2=3$ and
$dim_{R/\m} \m^2/\m^3=1$. Therefore, by \cite[Prop.2.8]{ama}, $\m$
can be generated by three elements, say, $\m=(x, y, z)$ for some
$x,y,z \in \m-\m^2$.

First we assume that $x\m=\{0\}$. Choose $w\in \m^2-\{0\}$ and
denote $u_1=x$, $u_2=w$, $u_3=x+w$, $v_1=y$, $v_2=z$, $v_3=y+z$,
$v_4=y+w$, $v_5=z+w$, $v_6=y+z+w$, and $v_7=x+y$; then $u_i\cdot
v_j=0$ for every $i, j$. It follows that $K_{3,7}\subseteq
\Gamma(R)$, so that $\gamma(\Gamma(R))\geq \gamma(K_{3,7})=2$.
Thus we may assume that $u\m\neq \{0\}$ for any $u\in \m-\m^2$.

Suppose that $x^2=y^2=z^2=0$. Since  $2ab=0$ for any $a,b \in \m$
since $\m^3=\{0\}$ and $|R/\m|=2$, we have that $u^2=0$ for all
$u\in\m$. Noting that $|\m^2|=2$, we may assume $xy\neq0$, that
is, $\m^2=\{xy,0\}$. Now if $xz\neq 0$, then $xz=xy$, so that
$x(z-y)=0$. Replacing $z$ by $z-y$, we may assume $xz=0$.
Furthermore, if $yz \neq 0$, then $yz=xy$, so that $y(z-x)=0$.
Since $x(z-x)=0$, we may replace $z$ by $z-x$ and assume that
$yz=0$. However, this implies $z\m=\{0\}$, which contradicts the
assumption that $u\m\neq \{0\}$ for any $u\in \m-\m^2$. Thus,
$u^2\neq 0$ for some $u\in \{x,y,z\}$.  After a suitable change of
$x,y,z$, we may assume that $x^2 \neq 0$ and $xy=xz=0$.  There are
thus two cases to consider, either
$y^2\neq0$ or $y^2=0$.\\
\textbf{Case 1} : $y^2\neq 0$. In this case, we may further assume
that $yz=0$ as in the previous paragraph. Therefore $z^2\neq 0$ as
$z\m\neq \{0\}$. Consequently, $x^2=y^2=z^2$ and $xy=xz=yz=0$. Let
$u_1=x^2$, $u_2=x$, $u_3=x+x^2$, $v_1=y$, $v_2=z$, $v_3=y+z$,
$v_4=y+x^2$, $v_5=z+x^2$, $v_6=y+z+x^2$, $w_1=x+y$, $w_2=x+y+x^2$,
$w_3=x+z$, $w_4=x+z+x^2$, $w_5=x+y+z$ and $w_6=x+y+z+x^2$.  Then
$u_i\cdot v_j=0$ for every $i, j$. It follows that
$K_{3,6}\subseteq \Gamma(R)$, and therefore
$\gamma(\Gamma(R))\geq\gamma(K_{3,6})= 1$. Let $G=\Gamma(R)$,
$G'=G-\{u_1w_1, u_1w_2,\dots,u_1w_6,u_1u_2, u_1u_3, v_1v_2,
v_1v_5, v_2v_4, v_3v_6\}$, and
$G''=G'-\{w_1,w_2,\dots,w_6\}$. We then have $G''\simeq K_{3,6}$. 

Suppose that $\gamma(G)=1$. Since
$1=\gamma(G'')\leq\gamma(G')\leq\gamma(G)$, we have
$\gamma(G')=1$. Since
$$V(G')=\{u_1, u_2, u_3, v_1, \dots, v_6, w_1, \dots, w_6\}$$ and
$$\ba{rl}E(G')=&\,\,\,\,\,\, \{u_iv_j|~1\leq i\leq 3, 1\leq j\leq 6\}
               \,\cup\, \{w_1w_2, w_3w_4\}
               \,\cup\, \{w_iw_j|~1\leq i\leq 4, j=5, 6\}\\
               &\cup\,\, \{v_1w_3, v_1w_4, v_2w_1, v_2w_2, v_3w_5, v_3w_6,
                       v_4w_3,v_4w_4, v_5w_1, v_5w_2, v_6w_5, v_6w_6 \},\ea$$
by the Euler characteristic formula, there are 25 faces when
drawing $G'$ on a torus. Fix a representation of $G'$ and let
$\{F'_1, \dots, F'_{25}\}$ be the set of faces of $G'$
corresponding to that representation. Since $G''\simeq K_{3,6}$,
this graph has 9 faces (see Remark \ref{rem1}). Let $F''_1, \dots,
F''_9$ be the faces of $G''$ obtained by deleting $w_1, \dots,
w_6$ and all edges incident with $w_1, \dots, w_6$ from this
representation. Again, $\{F'_1, \dots, F'_{25}\}$ can be recovered
by inserting $w_1, \dots, w_6$ and all edges incident with $w_1,
\dots, w_6$ into the representation corresponding to $\{F''_1,
\dots, F''_9\}$. We note that $w_iw_j\in E(G')$ for $i=1, \dots,
4$ and $j=5,6$.  Therefore all $w_i$ should be inserted into the
same face, say
$F''_m$, of $G''$ to avoid crossings. 
Moreover, since $v_1w_3,v_2w_4,v_3w_5,v_4w_4,v_5w_2,v_6w_5\in
E(G')$, it follows that $v_1,v_2,\dots,v_6$ are all vertices of
$F''_m$. This contradicts the fact that the boundary of $F''_m$ is
a 4-cycle.
Thus, we conclude that $\gamma(G)\geq 2$.\\

\noindent \textbf{Case 2} : $y^2=0$. In this case, we may further
assume that $z^2=0$. %
Since $z\m\neq \{0\}$ by assumption, we have $yz\neq 0$.
Consequently, $x^2=yz$ and $xy=xz=y^2=z^2=0$. Let $u_1=x^2$,
$u_2=x$, $u_3=x+x^2$, $v_1=y$, $v_2=z$, $v_3=y+z$, $v_4=y+x^2$,
$v_5=z+x^2$, $v_6=y+z+x^2$, $w_1=x+y$, $w_2=x+y+x^2$, $w_3=x+z$,
$w_4=x+z+x^2$, $w_5=x+y+z$ and $w_6=x+y+z+x^2$.  Then $u_i\cdot
v_j=0$ for every $i, j$, so that
$K_{3,6}\subseteq \Gamma(R)$. Therefore $\gamma(\Gamma(R))\geq\gamma(K_{3,6})= 1$. %
Let $G=\Gamma(R)$, $G'=G-\{w_5, w_6\}-\{u_1w_1, u_1w_2,u_1w_3,
u_1w_4, u_1u_2, u_1u_3, v_1v_4, v_2v_5, v_3v_6\}$, and
$G''=G'-\{w_1, w_2,w_3, w_4\}$.  Thus $G''\simeq K_{3,6}$.

Suppose that $\gamma(G)=1$, and so $\gamma(G')=1$. Since
$$V(G')=\{u_1, u_2, u_3, v_1, v_2,\dots, v_6, w_1, \dots, w_4\}$$ and
$$\ba{rl} E(G')= &\,\,\,\,\,\, \{u_iv_j~|~1\leq i\leq 3, 1\leq j\leq 6\}
                 \,\cup\, \{w_1w_3, w_1w_4, w_2w_3, w_2w_4\}\\
                 &\cup\,\, \{v_1w_1, v_1w_2,v_2w_3, v_2w_4, v_4w_1, v_4w_2, v_5w_3,
                 v_5w_4\},\ea$$
$G'$ has 17 faces. By a similar deletion and insertion argument as
before, fix a representation of $G'$ and let $\{F'_1, \dots,
F'_{17}\}$ be the set of faces of $G'$ corresponding to this
representation. Let $\{F''_1, \dots, F''_9\}$ be the set of faces
of $G''$ obtained by deleting $w_1, \dots, w_4$ and all edges
incident with $w_1, \dots, w_4$ from $G'$.  Since
$w_1w_3,w_2w_3,w_1w_4,w_2w_4\in E(G')$, $w_1,w_2,w_3,$ and $w_4$
should be inserted into the same face, say $F''_n$, of $G''$ in
the recovering process from $G''$ to $G'$ to avoid crossings.
We note that $v_1w_1, v_2w_3, v_4w_1, v_5w_4 \in E(G')$ and
therefore
$v_1,v_2,v_4,v_5$ are the four vertices of $F''_n$. %
Denote $e_1=v_1w_1$,$e_2=v_1w_2$, $e_3=v_2w_3$, $e_4=v_2w_4$,
$e_5=w_1w_3$, $e_6=w_1w_4$, $e_7=w_2w_3$, and $e_8=w_2w_4$. Then
we obtain Figure~2 by inserting $w_1, \dots, w_4$ and $e_1, \dots,
e_4$ into $F''_n$. However, from Figure~2 we see that there is no
way to insert $e_5, \dots, e_8$ into $F''_n$ without crossings, a
contradiction. Therefore, we conclude that $\gamma(G)\geq 2$.
\end{proof}
\begin{picture}(200,136)(60, -120)

\put(260,0){\line(1, -1){60}}\put(320,-60){\line(-1, -1){60}}
\put(260,-120){\line(-1,1){60}}\put(200,-60){\line(1,1){60}}
\put(257.36,-2.64){$\bullet$}\put(257.36,-122.64){$\bullet$}
\put(317,-62.32){$\bullet$}\put(198,-62.32){$\bullet$}

\put(240,-48){$\bullet$}\put(240,-74){$\bullet$}\put(274,-47.6){$\bullet$}\put(273,-76.6){$\bullet$}

\put(190,-65){$v_1$}\put(324,-65){$v_2$}
\put(263,3){$v_4$}\put(263,-125){$v_5$}

\put(240,-40){$w_1$}\put(240,-80){$w_2$}\put(274,-40){$w_3$}\put(274,-82){$w_4$}

\put(199,-60){\line(4,-1){42}}\put(199,-60){\line(3,1){42}}
\put(318,-60.6){\line(-3,-1){42}}\put(318,-59){\line(-3,1){42}}

\put(215,-50){$e_1$}\put(215,-72){$e_2$}
\put(298,-50){$e_3$}\put(298,-72){$e_4$}

\put(230,0){$F''_n:$}

\put(240,-136){Figure 2}
\end{picture}\\

\begin{prop}
\label{examlem3} Let $(R, \m)$ be a local ring. If $|R|=16$,
$|R/\m|=2$ and $\m^3=0$, then $\gamma(\Gamma(\mathbb{Z}_2\times
R))\geq 2$.
\end{prop}
\begin{proof}
By hypothesis, we have $|\m|=8$. If $|\m^2|=4$, then $\m$ is
principal, and so is $\m^2$. This implies that $|\m^3|=2$, a
contradiction. So $\m^2=0$ or $|\m^2|=2$. If $\m^2=0$, then the
seven non-zero elements of $\m$ are all zero-divisors of $R$.
Write $\m-\{0\}= \{a_1, \dots, a_7\}$ and let $u_i=(0, a_i)$ and
$v_i=(1, a_i)$ for $i=1, \dots, 7$. Then $u_i\cdot v_j=0$ for
every $1\leq i, j\leq 7$, so that $K_{7,7}\subseteq
\Gamma(\mathbb{Z}_2\times R)$, and it follows that
$\gamma(\Gamma(\mathbb{Z}_2\times R))\geq \gamma(K_{7,7})=7$.
Next, we assume that $|\m^2|=2$, so that $dim_{R/\m} \m/\m^2=2$
and $dim_{R/\m} \m^2/\m^3=1$. Let $\m=(x,y)$ for some $x,y\in
\m-\m^2$. We now consider two cases.

\noindent\textbf{Case 1}.  Suppose $x^2=y^2=0$. In this case,
$xy\neq 0$ as $\m^2\neq 0$. Let $u_1=(0, xy)$, $u_2=(1, xy)$,
$u_3=(1, 0)$, $v_1=(0, x)$, $v_2=(0, x+xy)$, $v_3=(0, y)$,
$v_4=(0, y+xy)$, $v_5=(0, x+y)$, $v_6=(0, x+y+xy)$, $w_1=(1, x)$,
$w_2=(1, x+xy)$, $w_3=(1, y)$, $w_4=(1, y+xy)$, $w_5=(1, x+y)$ and
$w_6=(1, x+y+xy)$; then $u_i\cdot v_j=0$ for every $i, j$, so that
$K_{3,6}\subseteq \Gamma(\mathbb{Z}_2\times R)$, it follows that
$\gamma(\Gamma(\mathbb{Z}_2\times R))\geq \gamma(K_{3,6})=1$.
Write $G=\Gamma(\mathbb{Z}_2\times R)$.  Note that if $z\notin \m$
then $(0,z) \in V(G)$ and $deg\,((0,z))=1$.  So $V(\widetilde{G})$
consists of all $u_i,v_j$ and $w_k$ as defined above. Let
$G'=\widetilde{G}-\{w_3, w_4,w_5,w_6\}-\{u_1u_2, u_1u_3, v_1v_2,
v_3v_4, v_5v_6\}$ and $G''=G'-\{w_1, w_2\}$, so that $G''\simeq
K_{3,6}$.

Suppose that $\gamma(G)=1$.  Since
$\gamma(G'')\leq\gamma(G')\leq\gamma(G)$, we have $\gamma(G')=1$.
Since
$$V(G')=\{u_1, u_2, u_3, v_1, \dots, v_6, w_1, w_2\}$$ and
$$
E(G')= \{u_iv_j~|~1\leq i\leq 3, 1\leq j\leq 6\}
       \,\cup\,\{v_iw_j~|~i, j=1, 2\}
       \,\cup\,\,\{u_1w_1, u_1w_2 \},$$
$G'$ has 13 faces. Fix a representation of $G'$ and let $\{F'_1,
\dots, F'_{13}\}$ be the set of faces of $G'$ corresponding to
this representation. Let $\{F''_1, \dots, F''_9\}$ be the set of
faces of $G''$ obtained by deleting $w_1, w_2$ and all edges
incident with $w_1, w_2$ from $G'$. Therefore there are faces
$F''_{t_1}, F''_{t_2}$ so that inserting $w_1,w_2$ and all edges
incident with $w_1,w_2$ into these faces, we are able to recover
the set of faces $\{F'_1, \dots, F'_{13}\}$. Since
$w_1v_1=w_2v_1=w_1v_2=w_2v_2\in V(G')$, $w_1$ and $w_2$ should be
inserted into the same face, say $F''_l$, of $G''$. Let
$e_1=v_1w_1$, $e_2=v_1w_2$, $e_3=v_2w_1$, $e_4=v_2w_2$,
$e_5=u_1w_1$ and $e_6=u_1w_2$. After inserting $w_1, w_2$ and
$e_1, e_2,e_3, e_4$ into $F'_l$ we obtain Figure 3 as below. From
the figure, we see that there is no way to insert $e_5, e_6$ into
$F'_l$ without any crossings. Thus, we
conclude that $\gamma(G)\geq 2$.\\

\begin{picture}(200,150)(-20,-150)
\put(160,0){\line(1, -1){60}}\put(220,-60){\line(-1, -1){60}}
\put(160,-120){\line(-1,1){60}}\put(100,-60){\line(1,1){60}}

\put(157.36,-2){$\bullet$} \put(157,-122){$\bullet$}
\put(158,-51){$\bullet$}\put(217,-63){$\bullet$}
\put(98,-63){$\bullet$} \put(157,-75){$\bullet$}
\put(97,-70){$v_1$}\put(220,-70){$v_2$} \put(170,0){$u_1$}

\put(220,-60.32){\line(-5, 1){60}}\put(220,-60.6){\line(-5, -1){60}}

\put(100,-60.32){\line(5, 1){62}}\put(100,-61){\line(5, -1){60}}
\put(130,-50){$e_1$}\put(130,-74){$e_2$}\put(180,-50){$e_3$}\put(180,-75){$e_4$}
\put(154,-41){$w_1$}\put(155,-80){$w_2$}

\put(120,0){$F'_l:$}

\put(136,-136){Figure~3}
\end{picture}\\
\noindent\textbf{Case 2}.  Suppose $x^2\neq 0$. In this case, we
may assume $xy=0$.  Otherwise, $x^2=xy$, and after replacing $y$
by $y-x$ we have $xy=0$. We note that either $y^2=0$ or $y^2=x^2$
as $|\m^2|=2$. Assume that $y^2=0$. Let $u_1=(0, y)$, $u_2=(0,
x^2)$, $u_3=(0, y+x^2)$, $v_1=(0, x+y)$, $v_2=(1, x+y)$, $v_3=(0,
x)$, $v_4=(1, x)$, $v_5=(0, x+x^2)$, $v_6=(0, x+y+x^2)$ and
$v_7=(1, x+x^2)$; then $u_i\cdot v_j=0$ for every $i, j$, so that
$K_{3,7}\subseteq \Gamma(\mathbb{Z}_2\times R)$, and it follows
that $\gamma(\Gamma(\mathbb{Z}_2\times R))\geq \gamma(K_{3,7})=2$.
Therefore, it remains to discuss the case when $y^2=x^2$.  Let
$u_1=(0, x^2)$, $u_2=(1, x^2)$, $u_3=(1, 0)$, $v_1=(0, y)$,
$v_2=(0, y+x^2)$, $v_3=(0, x)$, $v_4=(0, x+x^2)$, $v_5=(0, x+y)$,
$v_6=(0, x+y+x^2)$, $w_1=(1, x)$, $w_2=(1, x+x^2)$, $w_3=(1, y)$,
$w_4=(1, y+x^2)$, $w_5=(1, x+y)$ and $w_6=(1, x+y+x^2)$.  Then
$u_i\cdot v_j=0$ for every $i, j$, so that $K_{3,6}\subseteq
\Gamma(\mathbb{Z}_2\times R)$, and it follows that
$\gamma(\Gamma(\mathbb{Z}_2\times R))\geq \gamma(K_{3,6})=1$.
 Let $G=\Gamma(\mathbb{Z}_2\times R)$ and
$G^*=\widetilde{G}-\{w_3,w_4,w_5, w_6\}-\{u_1u_2, u_1u_3, v_1v_3,
v_1v_4, v_2v_3, v_2v_4, v_5v_6\}$. We observe that $G^*$ is
isomorphic to the graph $G'$ obtained in Case 1. Therefore
$\gamma(G^*)\geq 2$ and thus $\gamma(G)\geq 2$.
\end{proof}

\begin{prop}\label{Rpn} Let $(R,\m)$ be a local ring such that
$|R|=p^n$, where $p$ is prime and $n \in \mathbb{N}$. If
$\m^{n-1}\neq0$, then $\m$ is principal and $\Gamma(R)\simeq
\Gamma(\mathbb{Z}_{p^n})$. Moreover, for any ring $S$, one has that
$\Gamma(S\times R)\simeq \Gamma(S\times \mathbb{Z}_{p^n})$.
\end{prop}

\begin{proof}
We note that $\m^{k+1}\subsetneq \m^k$ if $\m^k \neq 0$, so that
$|\m^k|\leq p^{n-k}$ for $k\leq n$. Since $\m^{n-1}\neq 0$ we have
$|\m^k|=p^{n-k}$, so that $\dim_{R/\m}{\m/\m^2}=1$ and therefore
$\m$ is principal. Suppose $\m=(x)$ and let 
$\widetilde{\m}_k=\m^k-m^{k+1}$.  Then
$\widetilde{\m}_k=\{c_1x^k,\dots,c_lx^k\}$, where
$l=p^{n-k-1}(p-1)$ and each $c_i$ is a unit in $R$. Observe that
for $u_i\in \widetilde{\m}_i$ and $v_j\in \widetilde{\m}_j$,
$u_i\cdot v_j=0$
 if and only if $i+j\geq n$. %
Thus, $\Gamma(R)$ is uniquely determined. In particular,
$\Gamma(R)\simeq \Gamma(\mathbb{Z}_{p^n})$ as $\mathbb{Z}_{p^n}$
satisfies these assumptions. Moreover, let $\psi: R\to
\mathbb{Z}_{p^n}$ be a bijective map such that
$\psi(\widetilde{\m}_k)=(\widetilde{\n}_k)$ for each $k<n$, where
$\n=(p)$ is the maximal ideal of $\mathbb{Z}_{p^n}$. It is easy to
see that $\phi: S\times R \to S\times \mathbb{Z}_{p^n}$ with
$\phi(a,b)=(a,\psi(b)) $ induces an isomorphism between
$\Gamma(S\times R)$ and $\Gamma(S\times \mathbb{Z}_{p^n})$.
\end{proof}

We now briefly turn our attention to the case where $\Gamma(R)$ is
planar. In the following three examples, we show by explicit
representations that the zero-divisor graphs of the listed local
rings are planar.

\begin{example}\label{R27}
(a) $\Gamma(\mathbb{Z}_{27})$ is isomorphic to $G_1$ as in Figure
4-1. (b) $\Gamma(\mathbb{Z}_{16})$ is isomorphic to $G_2$ as shown
in Figure 4-2.
\end{example}

\begin{proof}
(a) In $\mathbb{Z}_{27}$, let $u_1=\bar 3$, $u_2=\bar 6$, $u_3=\bar
{12}$, $u_4=\bar {15}$, $u_5=\bar{21}$, $u_6=\bar{24}$,
$v_1=\bar{9}$, $v_2=\bar{18}$. Since $u_i\cdot v_j=0$  and $v_1\cdot
v_2=0$, $\Gamma(\mathbb{Z}_{27})$ is isomorphic to $G_1$ as shown in
Figure 4-1.

(b) In $\mathbb{Z}_{16}$, let $u=\bar 8$, $v_1=\bar 2$, $v_2=\bar
6$, $v_3=\bar {4}$, $v_4=\bar {12}$, $v_5=\bar{10}$,
$v_6=\bar{14}$. Since $u\cdot v_j=0$  and $v_3\cdot v_4=0$,
$\Gamma(\mathbb{Z}_{16})$ is isomorphic to $G_2$ as shown in
Figure 4-2.
\end{proof}

\begin{picture}(200,120)(0,-100)

\put(100,-15){$\bullet$}\put(100,-75){$\bullet$}
\put(100,-6){$v_1$}\put(100,-81){$v_2$}
\put(104,-11.6){\line(-5,-2){76}} \put(104,-11.6){\line(-4,-3){40}}
\put(103.6,-11.6){\line(-3,-5){18}} \put(102,-11.6){\line(3,-5){18}}
\put(102,-11.6){\line(4,-3){40}}  \put(102,-11.6){\line(5,-2){76}}

\put(102.8,-11.6){\line(0,-1){60}}

\put(104,-72.6){\line(-5,2){76}} \put(104,-72.6){\line(-4,3){40}}
\put(103.6,-72.6){\line(-3,5){18}} \put(102,-72.6){\line(3,5){18}}
\put(102,-72.6){\line(4,3){40}}  \put(102,-72.6){\line(5,2){76}}

\put(25,-45){$\bullet$} \put(61,-45){$\bullet$}
\put(83,-45){$\bullet$} \put(117.5,-45){$\bullet$}
\put(140,-45){$\bullet$} \put(176,-45){$\bullet$}

\put(15,-45){$u_1$} \put(51,-45){$u_2$} \put(75,-45){$u_3$}
\put(123,-45){$u_4$} \put(146,-45){$u_5$} \put(182,-45){$u_6$}
\put(70,-7){$G_1:$} \put(82,-95){Figure~4-1}

\put(300,-12){\line(1, -1){30}}\put(300,-12){\line(1, -2){15}}
\put(300,-12){\line(-1,-1){30}}
\put(300,-12){\line(2,-1){60}}\put(300,-12){\line(-2, -1){60}}

\put(300,-12){\line(-1, -2){15}}
\put(300,-7){$u$}\put(297.2,-15){$\bullet$} \put(228,-50){$v_1$}
\put(238,-45){$\bullet$}\put(268,-45){$\bullet$}
\put(282,-45){$\bullet$}\put(312.4,-45){$\bullet$}\put(284,-42){\line(1,0){29}}
\put(327,-45){$\bullet$}\put(357,-45){$\bullet$}

\put(265,-50){$v_2$}\put(280,-50){$v_3$}\put(310,-50){$v_4$}\put(325,-50){$v_5$}\put(355,-50){$v_6$}

\put(270,-7){$G_2 :$}\put(275,-92){Figure~4-2}

\end{picture}

\begin{example}\label{examlem4}
Let $R$ be one of the following local rings: (a)~${\mathbb{Z}_2[x,
y]}/{(x^3, xy, y^2-x^2)}$, (b)~${\mathbb{Z}_4[x]}/{(x^3,
x^2-2x)}$, (c)~${\mathbb{Z}_4[x,y]}/{(x^3, x^2-2, xy, y^2-2)}$, or
(d)~${\mathbb{Z}_8[x]}/{(x^2-4, 2x)}$.  Then $\Gamma(R)$ is planar
and is isomorphic to $G_3$ as shown in Figure~5-1.
\end{example}
\begin{proof} We briefly sketch the details for each case.

(a) If $R={\mathbb{Z}_2[x, y]}/{(x^3,xy, y^2-x^2)}$, then
$Z(R)^*=\{\bar x,\bar y,\bar x^2,\bar x+\bar x^2,\bar y+\bar x^2,
\bar x+\bar y,\bar x+\bar y+\bar x^2\}$. Let $u=\bar x^2$,
$v_1=\bar x$, $v_2=\bar y$, $v_3=\bar x+\bar x^2$, $v_4=\bar
y+\bar x^2$, $v_5=\bar x+\bar y$, and $v_6=\bar x+\bar y+\bar
x^2$. Then $\Gamma(R)\simeq G_3$ as shown in Figure~5-1.

(b) If $R={\mathbb{Z}_4[x]}/{(x^3, x^2-2x)}$, then $Z(R)^*=\{\bar
2,\bar x, \bar{2x}, \bar 2+ \bar{2x}, \bar{3x}, \bar 2+\bar x,
\bar 2+\bar{3x}\}$. Let $u= \bar{2x}$, $v_1=\bar x$, $v_2=\bar
2+\bar x$, $v_3=\bar{3x}$, $v_4=\bar 2+\bar{3x}$, $v_5=\bar 2$ and
$v_6=\bar 2+ \bar{2x}$.  It then follows that $\Gamma(R)\simeq
G_3$ as shown in Figure~5-1.

(c) If $R={\mathbb{Z}_4[x,y]}/{(x^3, x^2-2, xy, y^2-2)}$, then
$Z(R)^*=\{\bar 2,\bar x,\bar y,\bar x+\bar 2,\bar y+\bar 2,\bar
x+\bar y, \bar x+\bar y+\bar 2\}$. Let $u=\bar 2$, $v_1=\bar x$,
$v_2=\bar y$, $v_3=\bar x+\bar 2$, $v_4=\bar y+\bar 2$, $v_5=\bar
x+\bar y$ and $v_6=\bar x+\bar y+\bar 2$.  Then $\Gamma(R)\simeq
G_3$ as shown in Figure~5-1.

(d) If $R={\mathbb{Z}_8[x]}/{(x^2-4, 2x)}$, then $Z(R)^*=\{\bar 2,
\bar 4, \bar 6, \bar x, \bar 2+\bar x,\bar 4+\bar x,\bar 6+\bar
x\}$. Let
$u=\bar 4$, %
$v_1=\bar 2$, %
$v_2=\bar x$, %
$v_3=\bar 6$, %
$v_4=\bar 4+\bar x$, %
$v_5=\bar 2+\bar x$ %
and $v_6=\bar 6+\bar x$. Then $\Gamma(R)\simeq G_3$ as shown in
Figure~5-1.
\end{proof}

\begin{picture}(200,90)(-40,-90)
\put(70,-20){\line(3, -1){80}}\put(70,-20){\line(1, -1){28}}

\put(70,-20){\line(-3, -1){92.2}}
\put(70.36,-20){\line(-1,-1){30}}\put(70,-20){\line(-3, -2){45}}
\put(97,-48){\line(1,0){53}}\put(71,-81){\line(-1, 1){29.5}}
\put(71,-81){\line(-3, 1){92.2}}

\put(67.6,-23){$\bullet$} \put(67.6,-83.6){$\bullet$}
\put(95,-50.6){$\bullet$} \put(148.6,-50.6){$\bullet$}

\put(38,-52.84){$\bullet$}\put(22,-52.84){$\bullet$}\put(-23.4,-53){$\bullet$}
\put(-22,-50.4){\line(1,0){47}}\put(23,-50.4){\line(1,0){17}}
\put(70.2,-21){\line(0,-1){60}}
\put(70,-15){$u$}\put(-28,-58){$v_1$}\put(15,-58){$v_2$}\put(33,-58){$v_3$}
\put(68,-90){$v_4$}\put(90,-55){$v_5$}\put(145,-55){$v_6$}

\put(40,-15){$G_3 :$} \put(45,-100){Figure~5-1}

\put(270,-20){\line(1, -1){30}}\put(270,-20){\line(1, -2){15}}
\put(270,-20){\line(-1,-1){30}}
\put(270,-20){\line(2,-1){60}}\put(270,-20){\line(-2, -1){60}}

\put(270,-20){\line(-1, -2){15}}
\put(270,-15){$u$}\put(267.4,-23){$\bullet$} \put(198,-58){$v_1$}
\put(208,-53){$\bullet$}\put(237.4,-53){$\bullet$}\put(210,-50){\line(1,0){30}}
\put(252,-53){$\bullet$}\put(282.4,-53){$\bullet$}\put(254,-50){\line(1,0){29}}
\put(297.5,-53){$\bullet$}\put(327,-53){$\bullet$}\put(298,-50){\line(1,0){30}}

\put(235,-58){$v_2$}\put(250,-58){$v_3$}\put(280,-58){$v_4$}\put(295,-58){$v_5$}\put(325,-58){$v_6$}

\put(240,-15){$G_4 :$}\put(245,-100){Figure~5-2}

\end{picture}\\
\begin{example}\label{examlem4.5}
If $R$ is one of the following local rings: (a)~${\mathbb{Z}_4[x,
y]}/{(x^2, y^2, xy-2)}$, (b)~${\mathbb{Z}_2[x, y]}/{(x^2, y^2)}$,
or (c)~${\mathbb{Z}_4[x]}/{(x^2)}$, then $\Gamma(R)$ is planar and
is isomorphic to $G_4$ as shown in Figure~5-2.
\end{example}
\begin{proof}(sketch)

(a) If $R={\mathbb{Z}_4[x, y]}/{(x^2, y^2, xy-2)}$, then
$Z(R)^*=\{\bar 2, \bar x, \bar y, \bar x+\bar 2, \bar y+\bar 2,
\bar x+\bar y, \bar x+\bar y+\bar 2\}$. Let $u=\bar 2$, $v_1=\bar
x$, $v_2=\bar x+\bar 2$, $v_3=\bar y$, $v_4=\bar y+\bar 2$,
$v_5=\bar x+\bar y$ and $v_6=\bar x+\bar y+\bar 2$.  Then
$\Gamma(R)\simeq G_4$ as shown in Figure~5-2.

(b) If $R={\mathbb{Z}_2[x, y]}/{(x^2, y^2)}$, then $Z(R)^*=\{\bar
x, \bar y, \bar x\bar y, \bar x+\bar x\bar y, \bar y+\bar x\bar y,
\bar x+\bar y, \bar x+\bar y+\bar x\bar y\}$. Let $u=\bar x\bar
y$, $v_1=\bar x$, $v_2=\bar x+\bar x\bar y$, $v_3=\bar y$,
$v_4=\bar y+\bar x\bar y$, $v_5=\bar x+\bar y$ and $v_6=\bar
x+\bar y+\bar x\bar y$.  Thus $\Gamma(R)\simeq G_4$ as shown in
Figure~5-2.

(c) If $R={\mathbb{Z}_4[x]}/{(x^2)}$, then $Z(R)^*=\{\bar 2, \bar
x, \bar{2x}, \bar 2+\bar{2x}, \bar{3x}, \bar 2+\bar x, \bar
2+\bar{3x}\}$. Let $u=\bar{2x}$, $v_1=\bar 2$, $v_2=\bar 2+\bar
{2x}$, $v_3=\bar x$, $v_4=\bar{3x}$, $v_5=\bar 2+\bar x$ and
$v_6=\bar 2+\bar{3x}$.  Then $\Gamma(R)\simeq G_4$ as shown in
Figure~5-2.
\end{proof}

From Lemma~\ref{genuskn}, we see that complete graphs on 7
vertices are interesting in the genus one setting.  Thus, we
discuss some examples of finite local rings whose zero-divisor
graphs have seven vertices.

\begin{remark} \label{k34}
Suppose  $G$ is a simple graph such that $|V(G)|=7$ and $G$
contains a subgraph isomorphic to $K_{3,4}$, i.e.,
$K_{3,4}\subseteq G\subseteq K_7$. Then by Lemma \ref{genuskn},
\ref{genusbi} and Remark \ref {reduced}, we have
$1=\gamma(K_{3,4})\leq\gamma(G)\leq\gamma(K_7)=1$.  In particular,
for the 4-partite graph $K_{1,1,1,4}$ one has
$\gamma(K_{1,1,1,4})=1$.

(2) Consider $G_5=\Gamma(\mathbb{Z}_{32}$).  Then $V(G_5)=\{\bar
2, \bar 4,\bar 6,\dots,\bar {30}\}$. Let $u_1=\bar
8,u_2=\bar{16},u_3=\bar{24}$, $v_1=\bar{4}$, $v_2=\bar{12}$,
$v_3=\bar{20}$, $v_4=\bar{28}$, $w_1=\bar{2}$, $w_2=\bar{6}$,
$w_3=\bar{10}$, $w_4=\bar{14}$, $w_5=\bar{18}$, $w_6=\bar{22}$,
$w_7=\bar{26}$, and $w_8=\bar{30}$. We note that $deg(w_i)=1$ for
each $i$, so that $\widetilde{G_5}=G_5-\{w_1,\dots,w_8\}$.
Moreover, since $u_i\cdot u_j=0$, $u_i\cdot v_j=0$ and $v_i\cdot
v_j\neq 0$ for all pairs $i,j$, we have $\widetilde{G_5}\simeq
K_{1,1,1,4}$. Therefore, $\gamma(G_5)=\gamma(\widetilde{G_5})=1$
and thus $G_5$ is isomorphic to $K_{1,1,1,4}$ with 8 single
pendant edges.\footnote{ For the representation of the 4-partite
graph $K_{1,1,1,4}$, refer to Figure 5(b) in \cite{wa}.}
\end{remark}

\begin{example}\label{examlem5}
Let $R$ be one of the following local rings: %
(a) ${\mathbb{Z}_2[x,y]}/{(x^3, xy, y^2)}$, %
(b) ${\mathbb{Z}_4[x]}/{(x^3, 2x)}$, %
(c) ${\mathbb{Z}_4[x,y]}/{(x^3,x^2-2, xy, y^2)}$, or%
(d) ${\mathbb{Z}_8[x]}/{(x^2, 2x)}$.  Then $\Gamma(R)\simeq
K_{1,1,1,4}$ and $\gamma(\Gamma(R))=1$.
\end{example}

\begin{proof}
(a) If $R={\mathbb{Z}_2[x,y]}/{(x^3, xy, y^2)}$, then
$Z(R)^*=\{\bar x, \bar y, \bar x^2, \bar x+\bar x^2, \bar y+\bar
x^2, \bar x+\bar y, \bar x+\bar y+\bar x^2\}$. Let $u_1=\bar y$,
$u_2=\bar x^2$, $u_3=\bar y+\bar x^2$, $v_1=\bar x$, $v_2=\bar
x+\bar x^2$, $v_3=\bar x+\bar y$ and $v_4=\bar x+\bar y+\bar x^2$.
It is easy to see that for all $i,j$ we have $u_i\cdot u_j=0$,
$u_i\cdot v_j=0$, and $v_i\cdot v_j\neq 0$.  Thus $\Gamma(R)\simeq
K_{1,1,1,4}$ and $\gamma(\Gamma(R))=1$ by Remark~\ref{k34} (1).

(b) If $R={\mathbb{Z}_4[x]}/{(x^3, 2x)}$, then $Z(R)^*=\{2, \bar
x, 2+\bar x, \bar x^2, 2+\bar x^2, \bar x+\bar x^2, 2+\bar x+\bar
x^2\}$. Let $u_1=2$, $u_2=\bar x^2$, $u_3=2+\bar x^2$, $v_1=\bar
x$, $v_2=2+\bar x$, $v_3=\bar x+\bar x^2$ and $v_4=2+\bar x+\bar
x^2$. As in the previous case, it is clear that $\Gamma(R)\simeq
K_{1,1,1,4}$.

(c) If $R={\mathbb{Z}_4[x,y]}/{( x^3,x^2-2, xy, y^2)}$, then
$Z(R)^*=\{2, \bar x, \bar y, \bar x+2, \bar y+2, \bar x+\bar y,
\bar x+\bar y+2\}$. Let $u_1=2$, $u_2=\bar y$, $u_3=\bar y+2$,
$v_1=\bar x$, $v_2=\bar x+2$, $v_3=\bar x+\bar y$ and $v_4=\bar
x+\bar y+2$. Again, it follows immediately that $\Gamma(R)\simeq
K_{1,1,1,4}$.

(d) If $R={\mathbb{Z}_8[x]}/{(x^2, 2x)}$, then $Z(R)^*=\{2, 4, 6,
\bar x, 2+\bar x, 4+\bar x, 6+\bar x\}$. Let $u_1=4$, $u_2=\bar
x$, $u_3=4+\bar x$, $v_1=2$, $v_2=\bar x+2$, $v_3=6$ and $v_4=\bar
x+6$.  Once again, it is easy to see that $\Gamma(R)\simeq
K_{1,1,1,4}$ and $\gamma(\Gamma(R))=1$ by Remark~\ref{k34}.
\end{proof}

\medskip
We end this section by discussing how the genus number of
zero-divisor graphs behaves with respect to products.

\begin{example}\label{z2z3f4} Let
$R=\mathbb{Z}_2\times\mathbb{Z}_3\times\mathbb{F}_4$, and denote
the elements of $\mathbb{F}_4$ as $x_0, x_1, x_2,$ and $x_3$.  We
note that the reduction of $\Gamma(R)$, denoted by $G_6$, consists
of 11 vertices, namely $u_1=(0,0,x_1), u_2=(0,0,x_2),
u_3=(0,0,x_3), v_1=(1,0,x_0),
v_2=(1,1,x_0),v_3=(1,2,x_0),v_4=(0,1,x_0),v_5=(0,2,x_0),w_1=(1,0,x_1),w_2=(1,0,x_2),w_3=(1,0,x_3)$,
with edge set
$E(G_6)=\{\,u_iv_j\,|\,i=1,2,3;j=1,2,3,4,5\,\}\cup\{\,v_rv_j\,|\,r=4,5;
j=1,2,3\,\}\cup\{\,v_rw_j\,|\,r=4,5;j=1,2,3\,\}$. Since
$K_{3,5}\subseteq G_6$, we have that $\gamma(G_6)\geq 1$. The
following embedding explicitly shows that $\gamma(G_6)=1$.

\end{example}
\hspace{1.5in}
\includegraphics{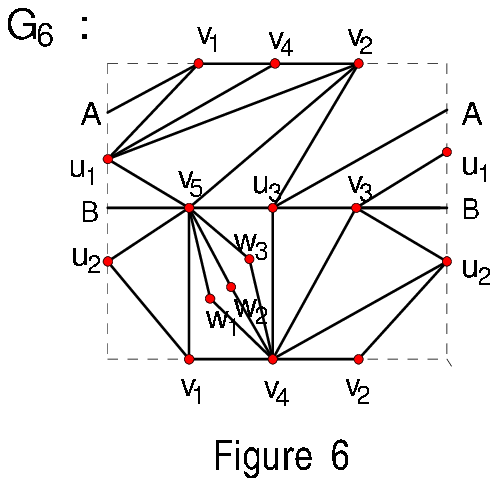}

\begin{prop}
\label{addexam} Let $(R, \m)$ be a finite local ring such that
$|\m|=4$ and $\m^2=0$.  Then $\gamma(\Gamma(\mathbb{Z}_2\times R))=1$. %
\end{prop}
\begin{proof}
Let $G=\Gamma(\mathbb{Z}_2\times R)$ and let $a_1, a_2, a_3$ be
distinct nonzero elements in $\m$. Let $u_i=(0, a_i)$ and $v_i=(1,
a_i)$ for $i=1, 2, 3$.  Then, $u_i\cdot v_j=0$ for each pair $i,
j$, so that $K_{3,3}\subseteq G$ and therefore $\gamma(G)\geq
\gamma(K_{3,3})=1$. On the other hand, let $w_1=(1, a_1)$,
$w_2=(0, a_1)$, $w_3=(1, a_2)$, $w_4=(0, a_2)$, $w_5=(1, a_3)$,
$w_6=(0, a_3)$, and $w_7=(1, 0)$.  Therefore
$V(\widetilde{G})=\{w_1, \dots, w_7\}$.
Let $\phi: \widetilde{G}\to \Gamma(\mathbb{Z}_{32})$ be the map
obtained by sending $w_i$ to $4i$ for every $i$. It is easy to see
that $\widetilde{G}$ is a subgraph of $\Gamma(\mathbb{Z}_{32})$
via $\phi$, so that $\gamma(G)=\gamma(\widetilde{G})\leq
\gamma(\Gamma(\mathbb{Z}_{32}))=1$ by \cite[Theorem~4.5]{wa}.
Thus, we conclude $\gamma(G)=1$.
\end{proof}

\begin{remark}\label{Fq} Let $\mathbb{F}_q$ denote the finite field with $q$
elements and let $\psi: \mathbb{F}_q \to \mathbb{Z}_q$ be any
bijective map such that $\psi(0)=0$. Let $R$ be a ring and define
a map $\phi : R\times \mathbb{F}_q \to R\times \mathbb{Z}_q$ such
that $\phi((a,b))=(a,\psi(b))$. Then $\phi$ induces an embedding
from $\Gamma(R\times \mathbb{F}_q)$ to $\Gamma(R\times
\mathbb{Z}_q)$. Therefore we conclude that $\gamma(\Gamma(R\times
\mathbb{F}_q))\leq \gamma(\Gamma(R\times \mathbb{Z}_q))$.
\end{remark}

\section{Toroidal zero-divisor graphs}
The main goal of this section is to determine all finite rings $R$
whose zero-divisor graphs are of genera at most one. To do this we
examine the characteristic of $R$, denoted by $char(R)$, the
cardinality of the residue field $|R/\m|$ (when $R$ is local), and
the number of irreducible components of $R$. We begin this section
with a few results related to the characteristic of a ring $R$.

\begin{lemma}
\label{keylem} Let $(R, \m)$ be a finite local ring such that
$|R|=q^2$ and $|R/\m|=q=p^d$ for some prime $p$ and
$d\in\mathbb{N}$. If $char(R)=p^2$, then $R\simeq
{\mathbb{Z}_{p^2}[x]}/{(f(x))}$, where $f$ is a monic polynomial of
$\mathbb{Z}_{p^2}[x]$ of degree $d$ such that the image of $f$ in
$\mathbb{Z}_p[x]$ is irreducible.
\end{lemma}
\begin{proof}
By hypothesis, we have $|\m|=q$. Since $\m/\m^2$ is a vector space
over $R/\m$, it follows that $|\m/\m^2|=q^r$, where
$r=\dim_{R/\m}\,\m/\m^2$. Therefore $\m^2=0$ and $r=1$, and this
implies that $\m=(x)$ for every $x\in\m-\{0\}$. Moreover, since
$char(R)=p^2$ and $\mathbb{Z}_{p^2}\subset R$, we have that $R$ is
a finitely generated $\mathbb{Z}_{p^2}$-algebra.  That is, $R=
\mathbb{Z}_{p^2}[\,u_1, \dots, u_n]$ for some $u_i \in R$.  We
note that $\m=(p)$ so that $\mathbb{F}_q \simeq R/\m =
\mathbb{Z}_{p^2}[\,u_1, \dots,
u_n]/(p)=\mathbb{Z}_p[\,\bar{u}_1,\dots,\bar{u}_n]$, where
$\bar{u}_i$ is the image of $u_i$ in $\mathbb{Z}_{p^2}[\,u_1,
\dots, u_n]/(p)$. Since $\mathbb{F}_q$ is a simple extension of
$\mathbb{Z}_p$, there exists $u\in R$ such that
$\mathbb{Z}_p[\,\bar{u}_1,\dots,\bar{u}_n]=\mathbb{Z}_p[\bar{u}]$.
Hence for any $w\in R$, there exists
$g_w(x)\in\mathbb{Z}_{p^2}[x]$
such that %
$$ w \equiv g_w(u)~~~~(mod ~\m).$$%
This implies that there exists $g_w^*(x_1,\dots,x_n) \in
\mathbb{Z}_{p^2}[\,x_1,\dots,x_n]$ such
that \be \label{eq1} w-g_w(u)-pg_w^*(u_1,\dots,u_n)=0. \ee %
Since $char(R)=p^2$, it follows from (\ref{eq1}) that%
\be \label{eq2} pw=pg_w(u). \ee %
That is, for every $w\in R$, there exists $g_w(x)\in
\mathbb{Z}_{p^2}[x]$ such that $pw=pg_w(u)$. Replacing $w$ by
$g_w^*(u_1,\dots,u_n)$, we see from (\ref{eq2}) that there exists
$g_w^{**}(x) \in \mathbb{Z}_{p^2}[x]$, such that
\be \label{eq3} pg_w^*(u_1,\dots,u_n)=pg_w^{**}(u).\ee %
It then follows from (\ref{eq1}) and (\ref{eq3}) that
$$w=g_w(u)+pg_w^{**}(u) \in \mathbb{Z}_{p^2}[u],$$
so that $R=\mathbb{Z}_{p^2}[u].$ Next, we proceed to show that there
exists a monic, d-degree polynomial $f(x) \in \mathbb{Z}_{p^2}[x]$
such that $f(u)=0$.

Since $[\,\mathbb{Z}_p[\bar u]:{\mathbb{Z}_p}]=d$, there exists a
monic, irreducible polynomial $h(x)\in\mathbb{Z}_p[x]$ of degree d
such that $h(\bar u)=0$. This implies that there exists a monic
polynomial \mbox{$l(x)\in\mathbb{Z}_{p^2}[x]$} of degree d such that
$\overline{l(u)}=h(\bar u)=0$, i.e., $l(u)\in \m$, so that
$l(u)=pl_1(u)$ for some $l_1(x)\in \mathbb{Z}_{p^2}[x]$. If
$\overline{l_1(u)}=0$ then $l_1(u)=pl_2(u)$ for some
$l_2(x)\in\mathbb{Z}_{p^2}[x]$. This implies that
$l(u)=p^2l_2(u)=0$. Set $f(x)=l(x)$ and then we are done. Otherwise,
since $[\,\mathbb{Z}_p[\bar u]:{\mathbb{Z}_p}]=d$,
$\overline{l_1(u)}=h_2(\bar u)$ for some $h_2(x)\in\mathbb{Z}_p[x]$
with $\deg h_2(x) < d$. Similarly, there exists $l_2(x)\in
\mathbb{Z}_{p^2}[x]$ such that $\deg l_2(x)=\deg h_2(x) < d$ and
$\overline{l_2(u)}=h_2(\bar u)$. Therefore
$\overline{l_1(u)}=\overline{l_2(u)}$ and this implies that
$l_1(u)-l_2(u)\in \m$, so that $pl_1(u)=pl_2(u)$. Let
$f(x)=l(x)-pl_2(x)\in\mathbb{Z}_{p^2}[x]$. Then
$f(u)=l(u)-pl_2(u)=l(u)-pl_1(u)=0$, and $\deg f(x)=d$ and the image
of $f$ in $\mathbb{Z}_p[x]$ (which is equal to $h(x)$) is
irreducible.

Let $\phi :\mathbb{Z}_{p^2}[x] \to \mathbb{Z}_{p^2}[u]$ be the
natural ring homomorphism which maps $x$ to $u$. We note that $\phi$
is surjective and that $f(x) \in \ker \phi$, so we have
$$q^2=|\hspace{0.8mm} \mathbb{Z}_{p^2}[u]\hspace{0.4mm}|=|\hspace{0.8mm}\mathbb{Z}_{p^2}[x]/\ker \phi \hspace{0.4mm}|\leq
|\hspace{0.8mm}\mathbb{Z}_{p^2}[x]/(f(x))\hspace{0.4mm}|.$$
However, since $\mathbb{Z}_{p^2}[x]/(f(x))$ is a free
$\mathbb{Z}_{p^2}$-module of rank d, it has $(p^2)^d=q^2$
elements. Thus $\ker \phi=(f(x))$ so that
$R=\mathbb{Z}_{p^2}[u]\simeq\mathbb{Z}_{p^2}[x]/(f(x))$ which
completes the proof.
\end{proof}

\medskip
Applying Lemma \ref{keylem} for the cases $(q,p,d)=(8,2,3)$ and
$(q,p,d)=(4,2,2)$ we obtain the following two corollaries.
\begin{coll}
\label{keylem1} Let $(R, \m)$ be a finite local ring such that
$|R|=64$ and $|R/\m|=8$. If $char(R)=4$, then $R\simeq
{\mathbb{Z}_4[x]}/{(x^3+x+1)}$. 
\end{coll}
\begin{proof}
We note that there are two irreducible polynomials in
$\mathbb{Z}_2[x]$ of degree 3, namely $x^3+x^2+1$ and $x^3+x+1$.
Applying Lemma~\ref{keylem} for $(q,p,d)=(4,2,2)$, we see that
$R\simeq \mathbb{Z}_4[x]/(f(x))$, where $\overline {f(x)}=x^3+x+1$
or $\overline {f(x)}=x^3+x^2+1$. Let $g_1(x)=x^3+x+1$,
$g_2(x)=x^3+x-1$, $g_3(x)=x^3+2x^2+x-1$, $g_4(x)=x^3+2x^2+x+1$,
$g_5(x)=x^3+x^2-1$, $g_6(x)=x^3+x^2+1$, $g_7(x)=x^3-x^2-1$,
$g_8(x)=x^3+x^2+1$, $h_1(x)=x^3-x+1$, $h_2(x)=x^3-x-1$,
$h_3(x)=x^3+2x^2-x-1$, $h_4(x)=x^3+2x^2-x+1$,
$h_5(x)=x^3-x^2+2x+1$, $h_6(x)=x^3-x^2+2x-1$,
$h_7(x)=x^3+x^2+2x+1$, and $h_8(x)=x^3+x^2+2x-1 \in
\mathbb{Z}_4[x]$. Then $f(x)$ is either $g_i[x]$ or $h_j[x]$. We
note that $\mathbb{Z}_4[x]/(g_i[x])$ are isomorphic by a linear
coordinate change and so are $\mathbb{Z}_4[x]/(h_j[x])$. Moreover,
consider the homomorphism $\sigma : \mathbb{Z}_4[x]/(x^3+x+1) \to
\mathbb{Z}_4[y]/(y^3-y+1)$ obtained by sending $\bar x$ to $\bar
y+2\bar y^2$. It is easy to see that $\sigma$ is an isomorphism,
and hence $R\simeq {\mathbb{Z}_4[x]}/{(x^3+x+1)}$.
\end{proof}

\begin{coll}
\label{keylem2} Let $(R, \m)$ be a finite local ring such that
$|R|=16$ and $|R/\m|=4$. If $char(R)=4$, then $R\simeq
{\mathbb{Z}_4[x]}/{(x^2+x+1)}$.
\end{coll}
\begin{proof}
We note that there is only one irreducible polynomial in
$\mathbb{Z}_2[x]$ of degree 2, namely $x^2+x+1$. Therefore,
applying Lemma~\ref{keylem} with $(q,p,d)=(4,2,2)$, we have that
$R\simeq \mathbb{Z}_4[x]/(f(x))$, where $\overline
{f(x)}=x^2+x+1$. Let $g_1(x)=x^2+x+1$, $g_2(x)=x^2+x+3$,
$g_3(x)=x^2+3x+1$ and $g_4(x)=x^2+3x+3 \in \mathbb{Z}_4[x]$. Then
$f(x)$ is one of these $g_i(x)$. Since $\mathbb{Z}_4[x]/(g_i(x))$
are isomorphic to each other, we conclude that $R\simeq
{\mathbb{Z}_4[x]}/{(x^2+x+1)}$.
\end{proof}

\medskip

We now recall a useful result from \cite{wa}.

\begin{theo}
\label{keytheo} Let $(R, \m)$ be a finite local ring which is not a
field. If $\gamma(\Gamma(R))\leq 1$, then $|R/\m|\leq 8$ and the
following hold:
\begin{description}
\item{(i)} If $|R/\m|= 8$, then $\m^2=0$ and $|R|=64$.
\item{(ii)} If $|R/\m|= 7$, then $\m^2=0$ and $|R|=49$.
\item{(iii)} If $|R/\m|= 5$, then $\m^2=0$ and $|R|=25$.
\item{(iv)} If $|R/\m|= 4$, then $\m^2=0$ and $|R|=16$.
\item{(v)} If $|R/\m|=3$, then $\m^3=0$ and $|R|\leq 27$.
\item{(vi)} If $|R/\m|=2$, then $\m^5=0$ and $|R|\leq 32$.
\end{description}
\end{theo}

\subsection*{Local rings with $\gamma(\Gamma(R))\leq 1$}

\begin{ntheo}\label{local}  From the above theorem we see that for a finite local ring
$(R,\m)$ with $\gamma(\Gamma(R))\leq 1$, the number of elements in
its residue field is bounded above by 8. Hence, we shall analyze
finite local rings by considering the cardinality of $R/\m$.
\end{ntheo}

\subsubsection*{$|R/\m|=8$.}

From Theorem~\ref{keytheo}, we have that $\m^2=0$ and $|R|=64$.
This implies that $|\m|=8$ and $\dim_{R/\m}\m/\m^2=1$ so that
$\m=(a)$ for every nonzero $a\in \m$ and that
$V(\Gamma(R))=\m-\{0\}$.  Therefore $\Gamma(R)\simeq K_7$ and then
$\gamma(\Gamma(R))=1$ by Lemma~\ref{genuskn}.  We note that if
$char(R)\geq 8$, then $\m=(2)$.  This implies that $char(R)=4$, a
contradiction. Thus, $char(R)\leq 4$.

\smallskip \noindent(i) If $char(R)=2$, then $R$ is equi-characteristic, so
that $R$ contains the field $\mathbb{F}_8$. Since $\m$ is principal,
it then follows that $R=\mathbb{F}_8[\m]\simeq
{\mathbb{F}_8[x]}/{(x^2)}$.

\medskip \noindent(ii) If $char(R)=4$, then from Corollary~\ref{keylem1} we
have that $R\simeq {\mathbb{Z}_4[x]}/{(x^3+x+1)}$. 

\subsubsection*{$|R/\m|=7$.}

From Theorem~\ref{keytheo}, we have that $\m^2=0$ and $|R|=49$.
This implies that $|\m|=7$ and $\dim_{R/\m}\m/\m^2=1$, so that
$\m=(a)$ for every nonzero $a\in \m$.  We note that since
$\m^2=0$, $\Gamma(R)\simeq K_6$, so that
$\gamma(\Gamma(R))=\gamma(K_6)=1$ by Lemma~\ref{genuskn}.   Since
$char(R)$ divides $|R|$, we have the following two cases to
consider.

\medskip \noindent(i) If $char(R)=7$, then $\mathbb{Z}_7\subseteq R$. Since
$\m$ is principal, $R=\mathbb{Z}_7[\m]\simeq
{\mathbb{Z}_7[x]}/{(x^2)}$.

\medskip \noindent(ii) If $char(R)=49$, then $\mathbb{Z}_{49}\subseteq R$. As
$|R|=49$, we conclude that $R=\mathbb{Z}_{49}$.

\subsubsection*{$|R/\m|=5$.}

From Theorem~\ref{keytheo}, we have that $\m^2=0$ and $|R|=25$.
This implies that $|\m|=5$ and $\dim_{R/\m}\m/\m^2=1$, so that
$\m=(a)$ for every nonzero $a\in \m$. Since $\m^2=0$, we have
$\Gamma(R)\simeq K_4$ and thus $\gamma(\Gamma(R))$ is planar by
Lemma~\ref{genuskn}.

\medskip \noindent(i) If $char(R)=5$, then  $R=\mathbb{Z}_5[\m]\simeq
{\mathbb{Z}_5[x]}/{(x^2)}$.

\medskip \noindent(ii) If $char(R)=25$, then $\mathbb{Z}_{25}\subseteq
R$ and thus $R=\mathbb{Z}_{25}$.

\subsubsection*{$|R/\m|=4$.}

From Theorem~\ref{keytheo}, we have that $\m^2=0$ and $|R|=16$. This
implies that $|\m|=4$ and $\dim_{R/\m}\m/\m^2=1$ so that $\m=(a)$
for every nonzero $a\in \m$. Since $\m^2=0$, we have
$\Gamma(R)\simeq K_3$, so that $\gamma(\Gamma(R))$ is planar. If
$char(R)\geq 8$, then $\m=(2)$, so that $char(R)=4$, a
contradiction. Thus $char(R)=2$ or $char(R)=4$.

\medskip \noindent(i) If $char(R)=2$, then $R$ is equi-characteristic, so
that $R$ contains the field $\mathbb{F}_4$. Since $\m$ is principal,
it follows that $R=\mathbb{F}_4[\m]\simeq
{\mathbb{F}_4[x]}/{(x^2)}$.

\medskip \noindent(ii) If $char(R)=4$, then $R\simeq
{\mathbb{Z}_4[x]}/{(x^2+x+1)}$ by Corollary~\ref{keylem2}.

\subsubsection*{$|R/\m|=3$.}

From Theorem~\ref{keytheo}, we have that $\m^3=0$ and $|R|\leq
27$, and therefore $|R|=9$ or $|R|=27$. We consider each case.

If $|R|=9$, then $\m^2=0$, so that $\Gamma(R)\simeq K_2$, which is
planar.  As in previous cases, we see that $R\simeq
\mathbb{F}_3[x]/(x^2)$ or $R\simeq \mathbb{Z}_9$.

If $|R|=27$, then $|\m|=9$. We note that $\m^2\neq 0$, for
otherwise $\Gamma(R)\simeq K_8$ which implies
$\gamma(\Gamma(R))=2$, a contradiction. It then follows from
Proposition~\ref{Rpn} and Example~\ref{R27} that $\Gamma(R)\simeq
\Gamma(\mathbb{Z}_{27})\simeq G_1$, as shown in Figure 4-1.

\medskip \noindent(i) If $char(R)=3$, then
$R=\mathbb{F}_3[\m]\simeq \mathbb{F}_3[x]/(x^3)$.

\medskip \noindent(ii) If $char(R)=9$, then $3 \in \m^2$.  Otherwise $\m=(3)$
and this implies that $\m^2=0$, a contradiction. Suppose that
$\m=(a)$ for some $a\in \m-\m^2$, so that either $3=a^2$ or
$3=-a^2$. Since $\m$ is principal, $R=\mathbb{Z}_9[a]$ and $R$ is
isomorphic to a quotient ring of either
$\mathbb{Z}_9[x]/(x^3,x^2-3)$ or $\mathbb{Z}_9[x]/(x^3,x^2+3)$.
However, both rings above have 27 elements. Thus we conclude that
$R\simeq\mathbb{Z}_9[x]/(x^3,x^2-3)$ or
$R\simeq\mathbb{Z}_9[x]/(x^3,x^2+3)$.

\medskip \noindent(iii) If $char(R)=27$, then $R=\mathbb{Z}_{27}$.

\subsubsection*{$|R/\m|=2$.}

From Theorem~\ref{keytheo}, we have that $\m^5=0$ and $|R|\leq
32$, and thus $|R|=4,8,16,$ or $32$.  As before, we consider cases
accordingly.

If $|R|=4$, then $\m^2=0$, so that $\Gamma(R)$ consists of a
single vertex. In this case, it is easy to see that $R\simeq
\mathbb{Z}_4$ or $R\simeq \mathbb{Z}_2[x]/(x^2)$.

Assume that $|R|=8$ and $\m^2=0$. Then $|\m|=4$ and $\dim_{R/\m}
\m/\m^2=2$.  Therefore $\Gamma(R)\simeq K_3$, which is a planar
graph.

\medskip \noindent(i) If $char(R)=2$, then  $R=\mathbb{Z}_2[\m]\simeq
\mathbb{Z}_2[x,y]/(x^2,xy,y^2)$.

\medskip \noindent(ii) If $char(R)=4$, then $\m=(2,a)$. We note that
$2a=a^2=0$.  Hence $R$ is isomorphic to a quotient ring of
$\mathbb{Z}_4[x]/(x^2,2x)$. Since there are 8 elements in
$\mathbb{Z}_4[x]/(x^2,2x)$, we conclude that
$R\simeq\mathbb{Z}_4[x]/(x^2,2x)$.

Assume that $|R|=8$ and $\m^2\neq 0$. By Proposition~\ref{Rpn}, we
see that $\Gamma(R)\simeq \Gamma(\mathbb{Z}_8)\simeq P_3$, the
path on 3 vertices.

\medskip \noindent(i) If $char(R)=2$, then $R$ contains the field
$\mathbb{Z}_2$, so that $R\simeq \mathbb{Z}_2[x]/(x^3)$.

\medskip \noindent(ii) If $char(R)=4$, then $\m^2=\{0,2\}$. Otherwise
$\m=(2)$ and this implies $\m^2=0$, which is a contradiction.
Therefore $\m=(a)$ with $a^2-2=a^3=0$. It follows that $R$ is
isomorphic to a quotient ring of $\mathbb{Z}_4[x]/(x^3,x^2-2)$.
However, $\mathbb{Z}_4[x]/(x^3,x^2-2)$ has 8 elements, so that
$R\simeq \mathbb{Z}_4[x]/(x^3,x^2-2)$.

\medskip \noindent(iii) If $char(R)=8$, then $R= \mathbb{Z}_8$.

Assume that $|R|=16$ and $\m^2=0$. It follows that $|\m|=8$ and
$\dim_{R/\m} \m/\m^2=3$.   Therefore, $\Gamma(R)\simeq K_7$ and
thus $\gamma(\Gamma(R))=\gamma(K_7)=1$ by Lemma~\ref{genuskn}.

\medskip \noindent(i) If $char(R)=2$, then  $R=\mathbb{Z}_2[\m]\simeq
\mathbb{Z}_2[x,y,z]/(x,y,z)^2$.

\medskip \noindent(ii) If $char(R)=4$, then $\m=(2,a,b)$, so that $R\simeq
\mathbb{Z}_4[x,y]/(x^2,xy,y^2,2x,2y)$.

Assume that $|R|=16$, $\m^2\neq 0$, and $\m^3=0$. If $|\m^2|=4$,
then $\m$ is principal and $|\m^3|=2$, which contradicts the
assumption $\m^3=0$. Therefore $|\m^2|=2$ and
$\dim_{R/\m}\m/\m^2=2$.

\medskip \noindent(i) Suppose $char(R)=2$.  Let $\m=(a,b)$ for some $a,b\in \m$. If
$a^2\neq 0$, then we may replace $b$ if necessary and assume that
$ab=0$ as $a^2$ is a basis element of $\m^2/\m^3$. We note that
$b^2=0$ or $b^2=a^2$.  Therefore $R$ is isomorphic to either
$\mathbb{Z}_2[x,y]/(x^3,xy,y^2)$ or
$\mathbb{Z}_2[x,y]/(x^3,xy,y^2-x^2)$. However, both rings in
question have 16 elements, so that either $R\simeq
\mathbb{Z}_2[x,y]/(x^3,xy,y^2)$ or $R\simeq
\mathbb{Z}_2[x,y]/(x^3,xy,y^2-x^2)$. In the first case,
$\Gamma(R)\simeq K_{1,1,1,4}$ and thus $\gamma(\Gamma(R))=1$ by
Example~\ref{examlem5}.  In the other case, $\Gamma(R)$ is planar
by Example~\ref{examlem4}. If $a^2=b^2=0$, then $R$ is isomorphic
to $\mathbb{Z}_2[x,y]/(x^2,y^2)$. Again, since
$\mathbb{Z}_2[x,y]/(x^2,y^2)$ has 16 elements, we conclude that
$R\simeq \mathbb{Z}_2[x,y]/(x^2,y^2)$, and thus $\Gamma(R)$ is
planar by Example~\ref{examlem4.5}.

\medskip \noindent(ii) Suppose $char(R)=4$.  First, we assume that $2\notin \m^2$.
Then $\m=(2,a)$ for some $a\in \m$. If $a^2\neq 0$, then $2a=0$ or
$2a=a^2$. By counting, we see that $R\simeq
\mathbb{Z}_4[x]/(x^3,2x)$ or $R\simeq
\mathbb{Z}_4[x]/(x^3,x^2-2x)$. In the first case,
$\gamma(\Gamma(R))=1$ by Example~\ref{examlem5} and in the second
case, $\Gamma(R)$ is planar by Example~\ref{examlem4}. If $a^2=0$,
then $R\simeq \mathbb{Z}_4[x]/(x^2)$, which is planar by
Example~\ref{examlem4.5}. Next, we assume that $2 \in \m^2$. Let
$\m=(a,b)$ for some $a,b \in \m$. From the remarks in (i), we see
that if $a^2\neq 0$, then we may assume that $ab=0$ which implies
that either $b^2=0$ or $b^2=2$. On the other hand, if $a^2=b^2=0$,
then $ab=2$ as $\m^2\neq 0$. By counting, we conclude that
$R\simeq \mathbb{Z}_4[x,y]/(x^3,x^2-2,xy,y^2)$ or $R \simeq
\mathbb{Z}_4[x,y]/(x^3,x^2-2,xy,y^2-2,y^3)$ or $R \simeq
\mathbb{Z}_4[x,y]/(x^2,xy-2,y^2)$. In the first case,
$\gamma(\Gamma(R))=1$ by Example~\ref{examlem5}.  In the other two
cases, we conclude $\Gamma(R)$ is planar by
Examples~\ref{examlem4} and ~\ref{examlem4.5} respectively.

\medskip \noindent(iii) Suppose $char(R)=8$. We note that $\m=(2,a)$ for some $a\in
\m$.  For otherwise $2\in \m^2$, which implies $4\in \m^4\neq 0$,
a contradiction. Since $4$ is the only nonzero element in $\m^2$,
we may replace $a$ by $a-2$ if necessary and assume that $2a=0$.
We note that $a^2=0$ or $a^2-4$, so that $R\simeq
\mathbb{Z}_8[x]/(x^2,2x)$ or $R\simeq\mathbb{Z}_8[x]/(x^2-4,2x)$.
In the first case, $\gamma(\Gamma(R))=1$ by Example~\ref{examlem5}
and $\Gamma(R)$ is planar in the second case by
Example~\ref{examlem4}.

Assume that $|R|=16$ and $\m^3\neq 0$. By Proposition~\ref{Rpn},
we see that $\m$ is principal and $\Gamma(R)\simeq
\Gamma(\mathbb{Z}_{16})$, which is planar by Example~\ref{R27}. As
before, we now consider cases depending on the characteristic of
$R$.

\medskip \noindent(i) Suppose $char(R)=2$. In this case, $R\simeq
\mathbb{Z}_2[x]/(x^4)$.

\medskip \noindent(ii) Suppose $char(R)=4$. We note that $2\in \m^2$, as otherwise
$\m=(2)$, and then $\m^2=0$, a contradiction.  Let $\m=(a)$ for
some $a\in \m$. Since $2a^2=0$, we see that $2=a^2$ or $2=a^3$ or
$2=a^2+a^3$. This implies that $R\simeq
\mathbb{Z}_4[x]/(x^2-2,x^4)$ or $R\simeq
\mathbb{Z}_4[x]/(x^3-2,x^4)$ or $R\simeq
\mathbb{Z}_4[x]/(x^3+x^2-2,x^4)$.

\medskip \noindent(iii) Suppose $char(R)=8$. We note that either $\m=(2)$ or $2\in
\m^2$, so that $\m^3=0$ or $4\in \m^4\neq 0$, both of which are
impossible.

\medskip \noindent(iv) Suppose $char(R)=16$. It is clear that $R=\mathbb{Z}_{16}$.

Assume that $|R|=32$. We note that $\m^2\neq 0$.  Otherwise
$\Gamma(R)\simeq K_{15}$. This implies that $\gamma(\Gamma(R))\geq
2$ by Lemma~\ref{genuskn}, which is a contradiction. Hence
$|\m^2|=2,4,$ or $8$. However, if $|\m^2|=2$ or $4$, then
$\gamma(\Gamma(R))\geq 2$ by Proposition~\ref{examlem1} and
Proposition~\ref{examlem2}. Therefore, $|\m^2|=8$, which implies
that $\m$ is principal and so is $\m^k$ for each $k\in
\mathbb{N}$. This implies that $\m^4\neq 0$. Thus,
$\Gamma(R)\simeq \Gamma(\mathbb{Z}_{32})$ by Proposition~\ref{Rpn}
and $\gamma(\Gamma(R))=1$ by \cite[Example~2.4]{wa}.

\medskip \noindent(i) Suppose $char(R)=2$. In this case, $R\simeq
\mathbb{Z}_2[x]/(x^5)$.

\medskip \noindent(ii) Suppose $char(R)=4$. We note that $2\in \m^3$, as otherwise
either $\m=(2)$ or $\m^2=(2)$ and both imply that $\m^4=0$, a
contradiction. Suppose that $\m=(a)$. Since $2a^2=0$, we see that
$2=a^3$ or $2=a^3+a^4=(a+a^2)^3$ or $2=a^4$. It follows that
$R\simeq \mathbb{Z}_4[x]/(x^3-2,x^5) \simeq
\mathbb{Z}_4[x]/((x+x^2)^3-2,x^5)$ or $R\simeq
\mathbb{Z}_4[x]/(x^4-2,x^5)$.

\medskip \noindent(iii) Suppose $char(R)=8$. If $\m=(2)$ or $2\in \m^3$, then
$\m^3=0$ or $4\in \m^6\neq0$, which are impossible. Thus,
$\m^2=(2)$. Let $\m=(a)$ for some $a\in \m$. Since $2a^2=4\in
\m^4$ and $4a=0\in \m^5$, we have that either $2=a^2$,
$2=-a^2=(a+2)^2$, $2=a^2+2a$, or $2=-a^2+2a$. Accordingly, it
follows that
$R\simeq\mathbb{Z}_8[x]/(x^2-2,x^5)\simeq\mathbb{Z}_8[x]/((x+2)^2-2,x^5)$
or $R\simeq \mathbb{Z}_8[x]/(x^2+2x-2,x^5)$ or $R\simeq
\mathbb{Z}_8[x]/(x^2-2x+2,x^5)$.

\medskip \noindent(iv) Suppose $char(R)=16$. We note either $\m=(2)$ or $2\in
\m^2$. This implies either $\m^4=0$ or $8\in\m^6\neq0$, both of
which are impossible.

\medskip \noindent(v) If $char(R)=32$, it is clear that $R=\mathbb{Z}_{32}$.

\medskip

Summarizing the above, we obtain the following theorems which
completely classify those local rings with a toroidal zero-divisor
graph. We are primarily concerned with the genus one case, but our
analysis recovers the following result from \cite{ns1}.  Note that
in that paper, some of the ideals are chosen with different
generators.  For instance, in that paper the rings
$\mathbb{Z}_9[x]/(x^2-3,3x)$ and $\mathbb{Z}_9[x]/(x^2-6,3x)$ are
listed as having planar zero-divisor graphs.  In our paper, we
use different generators for the ideals in question.

\begin{btheo}
\label{theo1} Let $(R, \m)$ be a finite local ring which is not a
field. Then $\Gamma(R)$ is planar if and only if $R$ is isomorphic
to one of the following 29 rings.

\medskip
\noindent
$\mathbb{Z}_4$, %
$\mathbb{Z}_8$, %
$\mathbb{Z}_9$, %
$\mathbb{Z}_{16}$, %
$\mathbb{Z}_{25}$, %
$\mathbb{Z}_{27}$, %
$\dfrac{\mathbb{Z}_2[x]}{(x^2)}$, %
$\dfrac{\mathbb{Z}_2[x]}{(x^3)}$, %
$\dfrac{\mathbb{Z}_2[x]}{(x^4)}$, %
$\dfrac{\mathbb{Z}_2[x,y]}{(x^2, y^2)}$
$\dfrac{\mathbb{Z}_2[x,y]}{(x^2,xy, y^2)}$, %
$\dfrac{\mathbb{Z}_2[x,y]}{(x^3,xy, y^2-x^2)}$, %
$\dfrac{\mathbb{F}_4[x]}{(x^2)}$, %
$\dfrac{\mathbb{Z}_3[x]}{(x^2)}$, %
$\dfrac{\mathbb{Z}_3[x]}{(x^3)}$, %
$\dfrac{\mathbb{Z}_4[x]}{(x^2)}$, %
$\dfrac{\mathbb{Z}_4[x]}{(x^2+x+1)}$, %
$\dfrac{\mathbb{Z}_4[x]}{(2x,x^2)}$, %
$\dfrac{\mathbb{Z}_4[x]}{(x^2-2,x^4)}$, %
$\dfrac{\mathbb{Z}_4[x]}{(x^3-2,x^4)}$, %
$\dfrac{\mathbb{Z}_4[x]}{(x^2-2,x^3)}$, %
$\dfrac{\mathbb{Z}_4[x]}{(x^3, x^2-2x)}$, %
$\dfrac{\mathbb{Z}_4[x]}{(x^3+x^2-2,x^4)}$, %
$\dfrac{\mathbb{Z}_4[x,y]}{(x^2, y^2, xy-2)}$, %
$\dfrac{\mathbb{Z}_4[x,y]}{(x^3,x^2-2, xy, y^2-2,y^3)}$, %
$\dfrac{\mathbb{Z}_5[x]}{(x^2)}$, %
$\dfrac{\mathbb{Z}_8[x]}{(x^2-4, 2x)}$, %
$\dfrac{\mathbb{Z}_9[x]}{(x^2-3,x^3)}$, %
$\dfrac{\mathbb{Z}_9[x]}{(x^2+3, x^3)}$. \\
\end{btheo}

\begin{btheo}\label{theo1-1} Let $(R, \m)$ be a finite local ring which is not a
field. Then $\gamma(\Gamma(R))=1$ if and only if $R$ is isomorphic
to one of the following 17 rings.

\medskip \noindent $\mathbb{Z}_{32}$, %
$\mathbb{Z}_{49}$, %
$\dfrac{\mathbb{Z}_2[x]}{(x^5)}$, %
$\dfrac{\mathbb{F}_8[x]}{(x^2)}$, %
$\dfrac{\mathbb{Z}_2[x,y]}{(x^3, xy, y^2)}$, %
$\dfrac{\mathbb{Z}_2[x, y, z]}{(x, y, z)^2}$, %
$\dfrac{\mathbb{Z}_4[x]}{(x^3+x+1)}$, %
$\dfrac{\mathbb{Z}_4[x]}{(x^3-2,x^5)}$, %
$\dfrac{\mathbb{Z}_4[x]}{(x^4-2, x^5)}$, %
$\dfrac{\mathbb{Z}_4[x,y]}{(x^3, x^2-2, xy, y^2)}$, %
$\dfrac{\mathbb{Z}_4[x]}{(x^3,2x)}$, %
$\dfrac{\mathbb{Z}_4[x,y]}{(2x ,2y, x^2,xy,y^2)}$, %
$\dfrac{\mathbb{Z}_7[x]}{(x^2)}$, %
$\dfrac{\mathbb{Z}_8[x]}{(x^2,2x)}$, %
$\dfrac{\mathbb{Z}_8[x]}{(x^2-2, x^5)}$, %
$\dfrac{\mathbb{Z}_8[x]}{(x^2+2x-2,x^5)}$, %
$\dfrac{\mathbb{Z}_8[x]}{(x^2-2x+2,x^5)}$.
\end{btheo}

\bigskip
\subsection*{Non-local rings with $\gamma(\Gamma(R))\leq 1$}

\begin{ntheo}\label{general}  Since a finite ring is Artinian,
it is isomorphic to a finite direct product of Artinian local
rings (see \cite[Theorem 8.7]{ama}). Thus, the number of maximal
ideals of $R$ is simply the number of components of $R$. From
\cite[Lemma~3.1]{wa}, we know that for a finite ring $R$,
$|Spec(R)|\geq 5$ implies that $\gamma(\Gamma(R))\geq 2$. Thus, to
find all rings $R$ with $\gamma(\Gamma(R))\leq 1$, we need only
consider all $R$ with $|Spec(R)|\leq 4$.
We thus consider cases according to the cardinality of $Spec(R)$.  As before, we will recover known results about the planar case as part of our analysis.%
\end{ntheo}

\subsubsection*{Let $|Spec(R)|=2$.}

By assumption, $R\simeq R_1\times R_2$ where both $(R_1, \m_1)$
and $(R_2, \m_2)$ are finite local rings with
$\gamma(\Gamma(R_i)\leq 1$. If $R_1$ and $R_2$ are both fields,
then $\Gamma(R)\simeq K_{m,n}$, where $m=|R_1|-1$ and $n=|R_2|-1$.
Therefore, by Lemma~\ref{genusbi}, $\Gamma(R)$ is planar if and
only if $R\simeq \mathbb{Z}_2\times \mathbb{F}_q$ or $R\simeq
\mathbb{Z}_3\times \mathbb{F}_q$, where $\mathbb{F}_q$ is a finite
field with $q$ elements.  Further, $\gamma(\Gamma(R))=1$ if and
only if $R\simeq \mathbb{F}_4\times \mathbb{F}_4$, $R\simeq
\mathbb{F}_4\times \mathbb{Z}_5$, $R\simeq \mathbb{F}_4\times
\mathbb{Z}_7$, or $R\simeq \mathbb{Z}_5\times \mathbb{Z}_5$. Thus,
we may assume that at least one of the $R_i$ is not a field. We
proceed by considering the pair $(|R_1|, |R_2|)$. Without loss, we
assume that $|R_1|\leq |R_2|$.

\medskip\noindent
\textbf{Case 1 : } $|R_1|=2$. In this case, $R_1\simeq
\mathbb{Z}_2$. Since $R_2$ is not a field, we have $|R_2/\m_2|\leq
8$ by Theorem~\ref{keytheo}, and also we see that if
$|R_2/\m_2|\geq 5$, then $\m_2^2=0$. Let $a_1, \dots, a_4$ be
distinct nonzero elements in $\m_2$ and let $u_i=(0, a_i)$,
$v_i=(1, a_i)$ for $i=1, \dots, 4$, and $v_5=(1, 0)$.  Then
$u_i\cdot v_j=0$ for every $i, j$, so that $K_{4,5}\subseteq
\Gamma(\mathbb{Z}_2\times R_2)$. It follows that
$\gamma(\Gamma(\mathbb{Z}_2\times R_2))\geq 2$.  We may thus
assume that $|R_2/\m_2|\leq 4$.  We now proceed by cases,
according to the cardinality
of $|R_2/\m_2|$.%

\medskip %
Assume that $|R_2/\m_2|=4$. Then we have that $|\m_2|=4$,
$\m_2^2=0$, and that $\m_2$ is principal by Theorem~\ref{keytheo}.
Therefore, $\gamma(\Gamma(\mathbb{Z}_2\times R_2))=1$ by
Proposition~\ref{addexam}. We see from the discussion in
(\ref{local}) that $R_2\simeq \dfrac{\mathbb{F}_4[x]}{(x^2)}$ or
$\dfrac{\mathbb{Z}_4[x]}{(x^2+x+1)}$, so that $R\simeq
\mathbb{Z}_2\times\dfrac{\mathbb{F}_4[x]}{(x^2)}$ or
$\mathbb{Z}_2\times\dfrac{\mathbb{Z}_4[x]}{(x^2+x+1)}$ in this
case.

\medskip %
Assume that $|R_2/\m_2|=3$. Then $\m_2^3=0$ and $|R_2|\leq 27$ by
Theorem~\ref{keytheo}.

\medskip\noindent
(i) If $|R_2|=3$, then $R_2\simeq \mathbb{Z}_3$ and
$\Gamma(\mathbb{Z}_2\times \mathbb{Z}_3)$ is planar by
\cite[Theorem~5.1]{afll}.

\medskip\noindent
(ii) If $|R_2|=9$, then $R_2\simeq \mathbb{Z}_9$ or
$R_2\simeq\dfrac{\mathbb{Z}_3[x]}{(x^2)}$. However,
$\Gamma(\mathbb{Z}_2\times \mathbb{Z}_9)\simeq
\Gamma(\mathbb{Z}_2\times \dfrac{\mathbb{Z}_3[x]}{(x^2)})$ by
Proposition~\ref{Rpn} and $\Gamma(\mathbb{Z}_2\times \mathbb{Z}_9)$
is planar by \cite[Theorem~5.1]{afll},  so $R\simeq
\mathbb{Z}_2\times \mathbb{Z}_9$ or $R\simeq \mathbb{Z}_2\times
\dfrac{\mathbb{Z}_3[x]}{(x^2)}$ in this case.

\medskip\noindent
(iii) If $|R_2|=27$, then $\m^2_2\neq 0$ from (\ref{local}), so
that $\Gamma(\mathbb{Z}_2\times R_2)\simeq
\Gamma(\mathbb{Z}_2\times \mathbb{Z}_{27})$ by
Proposition~\ref{Rpn}. However, $\gamma(\Gamma(\mathbb{Z}_2\times
\mathbb{Z}_{27}))\geq 2$ by \cite[Theorem~4.5]{wa}.

\medskip
Now, assume that $|R_2/\m_2|=2$. Then $\m_2^5=0$ and $|R_2|\leq
32$ by Theorem~\ref{keytheo}.

\medskip \noindent%
(i) If $|R_2|=4$, then $R_2\simeq \mathbb{Z}_4$ or $R_2\simeq
\dfrac{\mathbb{Z}_2[x]}{(x^2)}$. However, since
$\Gamma(\mathbb{Z}_2\times \mathbb{Z}_4)\simeq
\Gamma(\mathbb{Z}_2\times \dfrac{\mathbb{Z}_2[x]}{(x^2)})$ and
$\Gamma(\mathbb{Z}_2\times \mathbb{Z}_4)$ is planar by
\cite[Theorem~5.1]{afll}, we see that $R\simeq \mathbb{Z}_2\times
\mathbb{Z}_4$ or $R\simeq \mathbb{Z}_2\times
\dfrac{\mathbb{Z}_2[x]}{(x^2)}$.

\medskip \noindent%
(ii) Suppose that $|R_2|=8$ and $\m_2^2=0$. In this case,
$|\m_2|=4$, so that $\gamma(\Gamma(\mathbb{Z}_2\times R_2))=1$ by
Proposition~\ref{addexam}. Therefore, $R\simeq \mathbb{Z}_2\times
\dfrac{\mathbb{Z}_2[x, y]}{(x^2, xy, y^2)}$, or $R\simeq
\mathbb{Z}_2\times \dfrac{\mathbb{Z}_4[x]}{(2x, x^2)}$.

\medskip \noindent%
(iii) Suppose that $|R_2|=8$ and $\m_2^2\neq 0$. By
Proposition~\ref{Rpn} we have that $\Gamma(R_2)\simeq
\Gamma(\mathbb{Z}_8)$ and that $\Gamma(\mathbb{Z}_2\times
R_2)\simeq \Gamma(\mathbb{Z}_2\times \mathbb{Z}_8)$, which is
planar by \cite[Theorem~5.1]{afll}. Thus, $R\simeq
\mathbb{Z}_2\times \dfrac{\mathbb{Z}_2[x]}{(x^3)}$ or $R\simeq
\mathbb{Z}_2\times \dfrac{\mathbb{Z}_4[x]}{(x^2-2, x^3)}$ or
$R\simeq \mathbb{Z}_2\times \mathbb{Z}_8$ by (\ref{local}).

\medskip \noindent%
(iv) Suppose that $|R_2|=16$. In this case, $\m_2^3\neq 0$, as
otherwise $\gamma(\Gamma(R))\geq 2$ by Proposition~\ref{examlem3}.
Therefore, $\Gamma(R_2)\simeq \Gamma(\mathbb{Z}_{16})$ and
$\Gamma(R)\simeq \Gamma(\mathbb{Z}_2\times \mathbb{Z}_{16})$ by
Proposition~\ref{Rpn}. However, $\gamma(\Gamma(\mathbb{Z}_2\times
\mathbb{Z}_{16}))\geq 2$ by \cite[Example~2.5]{wa}.

\medskip \noindent%
(v) Finally, suppose that $|R_2|=32$. Then $\m_2^4\neq 0$ and
$\m_2^5=0$ from the discussion in (\ref{local}). Thus,
$\Gamma(R_2)\simeq \Gamma(\mathbb{Z}_{32})$ and $\Gamma(R)\simeq
\Gamma(\mathbb{Z}_2\times \mathbb{Z}_{32})$. However,
$\gamma(\Gamma(\mathbb{Z}_2\times \mathbb{Z}_{32}))\geq 2$ by
\cite[Example~2.5]{wa}.

\medskip \textbf{Case 2 : } $|R_1|=3$. In this case, $R_1\simeq
\mathbb{Z}_3$. Since $R_2$ is not a field, $|R_2/\m_2|\leq 8$ by
Theorem~\ref{keytheo}. If $|R_2/\m_2|\geq 4$, then $\m_2^2=0$ by
Theorem~\ref{keytheo}.  Let $a_1, a_2, a_3$ be distinct nonzero
elements in $\m_2$ and let $u_i=(0, a_i)$, $v_i=(1, a_i)$,
$v_{i+3}=(2, a_i)$ for $i=1, 2, 3$, and let $v_7=(1, 0)$.  Then
$u_i\cdot v_j=0$ for every $i, j$ so that $K_{3,7}\subseteq
\Gamma(\mathbb{Z}_3\times R_2)$ and therefore
$\gamma(\Gamma(\mathbb{Z}_3\times R_2))\geq 2$. 
We may thus assume henceforth that $|R_2/\m_2|\leq 3$.

\medskip %
Assume that $|R_2/\m_2|=3$. Then $\m_2^3=0$ and $|R_2|\leq 27$ by
Theorem~\ref{keytheo}.

\medskip \noindent%
(i) If $|R_2|=9$, then $R_2\simeq \mathbb{Z}_9$ or
$\dfrac{\mathbb{Z}_3[x]}{(x^2)}$. We note that
$\Gamma(\mathbb{Z}_3\times \mathbb{Z}_9)\simeq
\Gamma(\mathbb{Z}_3\times \dfrac{\mathbb{Z}_3[x]}{(x^2)})$ by
Proposition~\ref{Rpn} and $\Gamma(\mathbb{Z}_3\times
\mathbb{Z}_9)$ is planar by \cite[Theorem~5.1]{afll}.  Therefore,
$R\simeq \mathbb{Z}_3\times \mathbb{Z}_9$ or $R\simeq
\mathbb{Z}_3\times \dfrac{\mathbb{Z}_3[x]}{(x^2)}$ in this case.

\medskip \noindent%
(ii) If $|R_2|=27$, then $\m_2^2\neq 0$ by (\ref{local}), so that
$\Gamma(\mathbb{Z}_3\times R_2)\simeq \Gamma(\mathbb{Z}_3\times
\mathbb{Z}_{27})$ by Proposition~\ref{Rpn}. However, in this case
$\gamma(\Gamma(\mathbb{Z}_3\times \mathbb{Z}_{27}))\geq 2$.

\medskip %
Assume that $|R_2/\m_2|=2$. Then $\m_2^5=0$ and $|R_2|\leq 32$ by
Theorem~\ref{keytheo}.

\medskip \noindent%
(i) If $|R_2|=4$, then $R_2\simeq \mathbb{Z}_4$ or
$\dfrac{\mathbb{Z}_2[x]}{(x^2)}$. We note that
$\Gamma(\mathbb{Z}_3\times \mathbb{Z}_4)\simeq
\Gamma(\mathbb{Z}_3\times \dfrac{\mathbb{Z}_2[x]}{(x^2)})$ and
$\Gamma(\mathbb{Z}_3\times \mathbb{Z}_4)$ is planar by
\cite[Theorem~5.1]{afll}, so $R\simeq \mathbb{Z}_3\times
\mathbb{Z}_4$ or $R\simeq \mathbb{Z}_3\times
\dfrac{\mathbb{Z}_2[x]}{(x^2)}$ in this case.

\medskip \noindent%
(ii) Suppose that $|R_2|=8$ and $\m_2^2=0$. Let
$G=\Gamma(\mathbb{Z}_3\times R_2)$ and let $a_1, a_2, a_3$ be
distinct nonzero elements in $\m_2$. Let $u_i=(0, a_i)$, $v_i=(1,
a_i)$ and $v_{i+3}=(2, a_i)$ for $i=1, 2, 3$, and let $v_7=(1,
0)$. We then see that $u_i\cdot v_j=0$ for every $i, j$ as
$\m_2^2=0$ by Theorem~\ref{keytheo}, so that $K_{3,7}\subseteq
\Gamma(\mathbb{Z}_3\times R_2)$. It then follows that
$\gamma(\Gamma(\mathbb{Z}_3\times R_2))\geq 2$ by
Lemma~\ref{genusbi}.

\medskip \noindent%
(iii) Suppose that $|R_2|=8$ and $\m_2^2\neq 0$. Then $\Gamma(
R_2)\simeq \Gamma(\mathbb{Z}_8)$ and $\Gamma(\mathbb{Z}_3\times
R_2)\simeq \Gamma(\mathbb{Z}_3\times \mathbb{Z}_8)$ by
Proposition~\ref{Rpn}. We note that
$\gamma(\Gamma(\mathbb{Z}_3\times \mathbb{Z}_8))=1$ by
\cite[Theorem~4.5]{wa}. Thus either $R\simeq \mathbb{Z}_3\times
\dfrac{\mathbb{Z}_2[x]}{(x^3)}$, $R\simeq \mathbb{Z}_3\times
\dfrac{\mathbb{Z}_4[x]}{(x^2-2, x^3)}$, or $R\simeq
\mathbb{Z}_3\times \mathbb{Z}_8$ from the discussion in
(\ref{local}).

\medskip \noindent%
(iv) If $|R_2|\geq 16$, then $\Gamma(\mathbb{Z}_2\times R_2)$ is a
subgraph of $\Gamma(R)$. However, we see from Case~{1.C.}(iv) that
$\gamma(\Gamma(\mathbb{Z}_2\times R_2))\geq 2$.%

\medskip \textbf{Case 3 : } $|R_1|=4$.  In this case, $R_1\simeq
\mathbb{F}_4$ or $\mathbb{Z}_4$ or
$\dfrac{\mathbb{Z}_2[x]}{(x^2)}$. We note that if $|R_2|\geq 8$,
then $K_{3,7}\subseteq \Gamma(R)$, and so $\gamma(\Gamma(R))\geq
2$ by Lemma~\ref{genusbi}. Therefore, $|R_2|=4, 5,$ or $7$.
Equivalently, $R_2\simeq \mathbb{F}_4, \mathbb{Z}_4$,
$\dfrac{\mathbb{Z}_2[x]}{(x^2)}$, $\mathbb{Z}_5$, or
$\mathbb{Z}_7$. If $R_1\simeq \mathbb{F}_4$, then $R_2\simeq
\mathbb{Z}_4$ or $\dfrac{\mathbb{Z}_2[x]}{(x^2)}$ as $R_2$ is not
a field. It is easy to see that the reduction of
$\Gamma(\mathbb{F}_4\times \mathbb{Z}_4)$ is isomorphic to
$K_{3,3}$. Since $\gamma(K_{3,3})=1$ and
$\Gamma(\mathbb{F}_4\times \mathbb{Z}_4)\simeq
\Gamma(\mathbb{F}_4\times \dfrac{\mathbb{Z}_2[x]}{(x^2)})$, then
$R\simeq \mathbb{F}_4\times \mathbb{Z}_4$ or $\mathbb{F}_4\times
\dfrac{\mathbb{Z}_2[x]}{(x^2)}$. On the other hand, if $R_1\simeq
\mathbb{Z}_4$ and $R_2\simeq \mathbb{Z}_n$, then by
\cite[Theorem~3.5]{wa} we see that $\Gamma(\mathbb{Z}_4\times
\mathbb{Z}_n)$ is of genus one for $n=4, 5,$ and $7$. By
Proposition~\ref{Rpn}, we have $\Gamma(\mathbb{Z}_4\times S)\simeq
\Gamma(\dfrac{\mathbb{Z}_2[x]}{(x^2)}\times S)$ for any ring $S$.
We thus conclude that $R\simeq \mathbb{Z}_4\times \mathbb{Z}_4$,
$\mathbb{Z}_4\times\mathbb{Z}_5$, $\mathbb{Z}_4\times
\mathbb{Z}_7$, $\mathbb{Z}_4\times\dfrac{\mathbb{Z}_2[x]}{(x^2)}$,
$\mathbb{Z}_5\times\dfrac{\mathbb{Z}_2[x]}{(x^2)}$,
$\mathbb{Z}_7\times\dfrac{\mathbb{Z}_2[x]}{(x^2)}$, or
$\dfrac{\mathbb{Z}_2[x]}{(x^2)}\times\dfrac{\mathbb{Z}_2[x]}{(x^2)}$. %

\medskip
\textbf{Case 4 : } $|R_1|\geq 5$.  Then $|R_2|\geq 8$ as
$|R_2|\geq |R_1|$ by assumption. Therefore, $K_{4,7}\subseteq
\Gamma(R)$ and then $\gamma(\Gamma(R))\geq 2$ by
Lemma~\ref{genusbi}.

\subsubsection*{Let $|Spec(R)|=3$.}

Since $|Spec(R)|=3$, there are three finite local rings $(R_i,
\m_i)$, $i=1, 2, 3$ such that $R\simeq R_1\times R_2\times R_3$.
We proceed by considering the triple $(|R_1|, |R_2|, |R_3|)$.
Without loss we may assume that $|R_1|\leq |R_2|\leq |R_3|$, and
we consider the possibilities for such a triple.

\medskip
Suppose $|R_1|\geq 3$. In this case, $R_i\geq 3$ for each $i$. If
$|R_3|\geq 4$, then $K_{3,8}\subseteq \Gamma(R)$, which implies
that $\gamma(\Gamma(R))\geq 2$ by Lemma~\ref{genusbi}, a
contradiction. Therefore, we may assume that $|R_3|\leq 3$. Then
 $|R_i|=3$ for each $i$. Note that in this case,
$\gamma(\Gamma(\mathbb{Z}_3\times \mathbb{Z}_3\times
\mathbb{Z}_3))=1$ by \cite[Example~2.4]{wa}.

\medskip
Suppose $|R_1|=2$ and $|R_2|\geq 3$. In this case, if $|R_2|\geq
4$ or $|R_3|\geq 5$, then $K_{3,7}$ or $K_{4,5}$ is contained in
$\Gamma(R)$, so that $\gamma(\Gamma(R))\geq 2$ by
Lemma~\ref{genusbi}. Therefore, we may assume that $(|R_2|,
|R_3|)=(3, 3)$ or $(3, 4)$. If $(|R_2|, |R_3|)=(3, 3)$, then
$R\simeq \mathbb{Z}_2\times \mathbb{Z}_3\times \mathbb{Z}_3$. In
this case we have $\gamma(\Gamma(\mathbb{Z}_2\times
\mathbb{Z}_3\times \mathbb{Z}_3))=1$ by \cite[Theorem~5.1]{afll}
and $\Gamma(\mathbb{Z}_2\times \mathbb{Z}_3\times \mathbb{Z}_3)$
is a subgraph of $\Gamma(\mathbb{Z}_3\times \mathbb{Z}_3\times
\mathbb{Z}_3)$. Now, if $(|R_2|, |R_3|)=(3, 4)$ then
$R_2=\mathbb{Z}_3$ and $R_3=\mathbb{F}_4$, $\mathbb{Z}_4$, or
$\dfrac{\mathbb{Z}_2[x]}{(x^2)}$. %
Example \ref{z2z3f4} shows that $\gamma(\Gamma(\mathbb{Z}_2\times
\mathbb{Z}_3\times \mathbb{F}_4))=1$. However, by \cite[Example
2.5]{wa}, we have $\gamma(\Gamma(\mathbb{Z}_2\times
\mathbb{Z}_3\times \mathbb{Z}_4)) \geq 2$ and
$\gamma(\Gamma(\mathbb{Z}_2\times \mathbb{Z}_3\times
\dfrac{\mathbb{Z}_2[x]}{(x^2)})) \geq 2$ as well, since
$\Gamma(\mathbb{Z}_2\times \mathbb{Z}_3\times
\dfrac{\mathbb{Z}_2[x]}{(x^2)})\simeq \Gamma(\mathbb{Z}_2\times
\mathbb{Z}_3\times \mathbb{Z}_4)$ by Proposition~\ref{Rpn}.
Therefore we must have $R\simeq \mathbb{Z}_2\times
\mathbb{Z}_3\times \mathbb{Z}_3$ or $R\simeq\mathbb{Z}_2\times
\mathbb{Z}_3\times \mathbb{F}_4$.

\medskip
Suppose $|R_1|=2$ and $|R_2|=2$. In this case, if $|R_3|\geq 8$,
then $K_{3,7}$ is contained in $\Gamma(R)$ which forces
$\gamma(\Gamma(R))\geq 2$ by Lemma~\ref{genusbi}. Therefore, we
may assume that $|R_3|\leq 7$.  Note that $\mathbb{Z}_2\times
\mathbb{Z}_2\times \mathbb{Z}_n$ is planar when $n=2$ or $3$ by
\cite[Theorem~5.1]{afll} and $\gamma(\Gamma(\mathbb{Z}_2\times
\mathbb{Z}_2\times \mathbb{Z}_n))=1$ when $n=5$ or $7$ by
\cite[Theorem~3.5]{wa}.  Thus, it remains only to examine the case
where $|R_3|=4$, that is where $R_3=\mathbb{F}_4$, $\mathbb{Z}_4$,
or $ \dfrac{\mathbb{Z}_2[x]}{(x^2)}$. We note that
$\Gamma(\mathbb{Z}_2\times \mathbb{Z}_2\times
\dfrac{\mathbb{Z}_2[x]}{(x^2)})\simeq \Gamma(\mathbb{Z}_2\times
\mathbb{Z}_2\times \mathbb{Z}_4)$ by Proposition~\ref{Rpn} and
$K_{3,3}\subseteq \Gamma(\mathbb{Z}_2\times \mathbb{Z}_2\times
\mathbb{F}_4)\subseteq \Gamma(\mathbb{Z}_2\times
\mathbb{Z}_2\times \mathbb{Z}_4)$ by Remark~\ref{Fq}. Since
$\gamma(\Gamma(\mathbb{Z}_2\times \mathbb{Z}_2\times
\mathbb{Z}_4))=1$ by \cite[Example~2.4]{wa}, we see that
$\gamma(\Gamma(\mathbb{Z}_2\times \mathbb{Z}_2\times
\mathbb{Z}_4))=\gamma(\Gamma(\mathbb{Z}_2\times \mathbb{Z}_2\times
\dfrac{\mathbb{Z}_2[x]}{(x^2)})=\gamma(\Gamma(\mathbb{Z}_2\times
\mathbb{Z}_2\times \mathbb{F}_4))=1$.
Therefore, %
$R\simeq \mathbb{Z}_2\times \mathbb{Z}_2\times \mathbb{F}_4$ or
$\mathbb{Z}_2\times \mathbb{Z}_2\times \mathbb{Z}_4$ or $
\mathbb{Z}_2\times \mathbb{Z}_2\times
\dfrac{\mathbb{Z}_2[x]}{(x^2)}$ in this case.

\subsubsection*{Suppose $|Spec(R)|=4$.}

We assume that $R\simeq R_1\times R_2\times R_3\times R_4$ where
each $(R_i, \m_i)$ is a finite local ring for $i=1, 2, 3, 4$.  As
before we also assume that $|R_1|\leq |R_2|\leq |R_3|\leq|R_4|$.

\medskip
\textbf{Case 1 : } $|R_3|\geq 3$. In this case, certainly
$|R_4|\geq 3$. Therefore, $K_{3,8}\subseteq \Gamma(R)$, and then
$\gamma(\Gamma(R))\geq 2$ by Lemma~\ref{genusbi}.

\medskip
\textbf{Case 2 : } $|R_1|=|R_2|=|R_3|=2$. In this case, if
$|R_4|\geq 4$, then $K_{3,7}\subseteq \Gamma(R)$, so that
$\gamma(\Gamma(R))\geq 2$ by Lemma~\ref{genusbi}, a contradiction.
Therefore, $|R_4|\leq 3$, so that $R\simeq \mathbb{Z}_2\times
\mathbb{Z}_2\times \mathbb{Z}_2\times \mathbb{Z}_2$ or $R\simeq
\mathbb{Z}_2\times \mathbb{Z}_2\times \mathbb{Z}_2\times
\mathbb{Z}_3$. However, $\gamma(\Gamma(\mathbb{Z}_2\times
\mathbb{Z}_2\times \mathbb{Z}_2\times \mathbb{Z}_2))=1$ and
$\gamma(\Gamma(\mathbb{Z}_2\times \mathbb{Z}_2\times
\mathbb{Z}_2\times \mathbb{Z}_3))\geq 2$ by \cite[Example~2.5]{wa},
so that $R\simeq \mathbb{Z}_2\times \mathbb{Z}_2\times
\mathbb{Z}_2\times \mathbb{Z}_2$ in this case.

\medskip
We summarize the above work to obtain two theorems. These two
results provide a complete list of non-local rings whose
zero-divisor graphs have genera at most one.  The planar case has
been discussed previously in \cite{amy}, \cite{rb} and in
\cite{ns1}. Also, the next theorem appears in \cite{ns1}; we state
the result since it follows immediately from our analysis.  As
before, note that in \cite{ns1} different generators are chosen for
some of the ideals.

\bigskip
\begin{btheo} \label{theo2} Let $R$ be a
finite ring which is not local; then $\Gamma(R)$ is planar if and
only if $R$ is isomorphic to one of the following 15 types of rings.

\noindent
$\mathbb{Z}_2\times \mathbb{F}_q$, 
$\mathbb{Z}_3\times \mathbb{F}_q$, %
$\mathbb{Z}_2\times \mathbb{Z}_9$, %
$\mathbb{Z}_2\times \dfrac{\mathbb{Z}_3[x]}{(x^2)}$, %
$\mathbb{Z}_2\times \mathbb{Z}_4$, %
$\mathbb{Z}_2\times \dfrac{\mathbb{Z}_2[x]}{(x^2)}$, %
$\mathbb{Z}_2\times \dfrac{\mathbb{Z}_2[x]}{(x^3)}$, %
$\mathbb{Z}_2\times \dfrac{\mathbb{Z}_4[x]}{(x^2-2, x^3)}$, %
$\mathbb{Z}_2\times \mathbb{Z}_8$, %
$\mathbb{Z}_3\times \mathbb{Z}_9$, %
$\mathbb{Z}_3\times \dfrac{\mathbb{Z}_3[x]}{(x^2)}$, %
$\mathbb{Z}_3\times \mathbb{Z}_4$, %
$\mathbb{Z}_3\times \dfrac{\mathbb{Z}_2[x]}{(x^2)}$, %
$\mathbb{Z}_2\times \mathbb{Z}_2\times \mathbb{Z}_2$,
$\mathbb{Z}_2\times \mathbb{Z}_2\times \mathbb{Z}_3$.
\end{btheo}

\begin{btheo} \label{theo2-1}Let $R$ be a
finite ring which is not local; then $\gamma(\Gamma(R))=1$  if and
only if $R$ is isomorphic to one of the following 29 rings.

\medskip\noindent
$\mathbb{F}_4\times \mathbb{F}_4$, %
$\mathbb{F}_4\times \mathbb{Z}_5$, %
$\mathbb{F}_4\times \mathbb{Z}_7$, %
$\mathbb{Z}_5\times \mathbb{Z}_5$, %
$\mathbb{Z}_2\times \dfrac{\mathbb{F}_4[x]}{(x^2)}$, %
$\mathbb{Z}_2\times \dfrac{\mathbb{Z}_4[x]}{(x^2+x+1)}$, %
$\mathbb{Z}_2\times \dfrac{\mathbb{Z}_2[x, y]}{(x^2, xy, y^2)}$, %
$\mathbb{Z}_2\times \dfrac{\mathbb{Z}_4[x]}{(2x, x^2)}$, %
$\mathbb{Z}_3\times \dfrac{\mathbb{Z}_2[x]}{(x^3)}$, %
$\mathbb{Z}_3\times \dfrac{\mathbb{Z}_4[x]}{(x^2-2, x^3)}$, %
$\mathbb{Z}_3\times \mathbb{Z}_8$, %
$\mathbb{Z}_4\times \mathbb{F}_4$, %
$\mathbb{F}_4\times \dfrac{\mathbb{Z}_2[x]}{(x^2)}$, %
$\mathbb{Z}_4\times \mathbb{Z}_4$, %
$\mathbb{Z}_4\times \dfrac{\mathbb{Z}_2[x]}{(x^2)}$, %
$\dfrac{\mathbb{Z}_2[x]}{(x^2)}\times \dfrac{\mathbb{Z}_2[x]}{(x^2)}$, %
$\mathbb{Z}_4\times \mathbb{Z}_5$, %
$\mathbb{Z}_5\times \dfrac{\mathbb{Z}_2[x]}{(x^2)}$, %
$\mathbb{Z}_4\times \mathbb{Z}_7$, %
$\mathbb{Z}_7\times \dfrac{\mathbb{Z}_2[x]}{(x^2)}$, %
$\mathbb{Z}_2\times \mathbb{Z}_3\times \mathbb{Z}_3$, %
$\mathbb{Z}_3\times \mathbb{Z}_3\times \mathbb{Z}_3$, %
$\mathbb{Z}_2\times \mathbb{Z}_3\times \mathbb{F}_4$, %
$\mathbb{Z}_2\times \mathbb{Z}_2\times \mathbb{Z}_5$, %
$\mathbb{Z}_2\times \mathbb{Z}_2\times \mathbb{Z}_7$, %
$\mathbb{Z}_2\times \mathbb{Z}_2\times \mathbb{F}_4$, %
$\mathbb{Z}_2\times \mathbb{Z}_2\times \mathbb{Z}_4$, %
$\mathbb{Z}_2\times \mathbb{Z}_2\times \dfrac{\mathbb{Z}_2[x]}{(x^2)}$, and %
$\mathbb{Z}_2\times \mathbb{Z}_2\times \mathbb{Z}_2\times \mathbb{Z}_2$. %
\end{btheo}

\bigskip

At the end of this paper, we summarize all the results obtained
through the discussion in (\ref{local}), (\ref{general}) in four
tables showed in the following pages.

\section*{Acknowledgement}
\indent The authors wish to express their deepest gratitude to the
referee for careful reading of the paper, helpful comments and many
valuable suggestions.

\pagebreak

\tiny \centering
\begin{sideways}
\footnotesize $G_1,G_2,G_3,G_4$ are graphs as shown in Figure 4-1, 4-2, 5-1, 5-2. %
\end{sideways}
\begin{sideways}
\begin{tabular}{|c|cccccccccc|}
\hline
&&&&&&&&&&\\
$R$& %
$\mathbb{Z}_4$& %
$\dfrac{\mathbb{Z}_2[x]}{(x^2)}$& %

$\,\,\mathbb{Z}_9\,\,$& %
$\dfrac{\mathbb{Z}_3[x]}{(x^2)}$& %
$\mathbb{Z}_8$& %
$\dfrac{\mathbb{Z}_2[x]}{(x^3)}$& %
$\dfrac{\mathbb{Z}_4[x]}{(x^2-2,x^3)}$& %

$\dfrac{\mathbb{Z}_2[x,y]}{(x^2,xy, y^2)}$& %
$\dfrac{\mathbb{Z}_4[x]}{(2x,x^2)}$& %
$\dfrac{\mathbb{F}_4[x]}{(x^2)}$\\%

&&&&&&&&&&\\
\hline
&&&&&&&&&&\\
$|R|$&4&4&9&9&8&8&8&8&8&16\\
&&&&&&&&&&\\
$|R/\m|$&2&2&3&3&2&2&2&2&2&4\\
&&&&&&&&&&\\
$char(R)$&4&2&9&3&8&2&4&2&4&2\\
&&&&&&&&&&\\
$\Gamma(R)$&Point&Point&$K_2$&$K_2$&$P_3$&$P_3$&$P_3$&$K_3$&$K_3$&$K_3$\\
&&&&&&&&&&\\
$|V(\Gamma(R))|$&1&1&2&2&3&3&3&3&3&3\\
&&&&&&&&&&\\
$\gamma(\Gamma(R))$&0&0&0&0&0&0&0&0&0&0\\
&&&&&&&&&&\\
\hline
&&&&&&&&&&\\
$R$& %
$\dfrac{\mathbb{Z}_4[x]}{(x^2+x+1)}$& %
$\mathbb{Z}_{25}$& %
$\dfrac{\mathbb{Z}_5[x]}{(x^2)}$& %
$\mathbb{Z}_{16}$& %
$\dfrac{\mathbb{Z}_2[x]}{(x^4)}$& %
$\dfrac{\mathbb{Z}_4[x]}{(x^2-2,x^4)}$& %
$\dfrac{\mathbb{Z}_4[x]}{(x^3-2,x^4)}$& %
$\dfrac{\mathbb{Z}_4[x]}{(x^3+x^2-2,x^4)}$& %

$\dfrac{\mathbb{Z}_2[x,y]}{(x^3,xy,y^2-x^2)}$& %
$\dfrac{\mathbb{Z}_4[x]}{(x^3,x^2-2x)}$\\ %
&&&&&&&&&&\\
\hline
&&&&&&&&&&\\
$|R|$&16&25&25&16&16&16&16&16&16&16\\
&&&&&&&&&&\\
$|R/\m|$&4&5&5&2&2&2&2&2&2&2\\
&&&&&&&&&&\\
$char(R)$&4&25&5&16&2&4&4&4&2&4\\
&&&&&&&&&&\\
$\Gamma(R)$&$K_3$&$K_4$&$K_4$&$G_2$&$G_2$&$G_2$&$G_2$&$G_2$&$G_3$&$G_3$\\
&&&&&&&&&&\\
$|V(\Gamma(R))|$&3&4&4&7&7&7&7&7&7&7\\
&&&&&&&&&&\\
$\gamma(\Gamma(R))$&0&0&0&0&0&0&0&0&0&0\\
&&&&&&&&&&\\
\hline
&&&&&&&&&&\\
$R$& %

$\dfrac{\mathbb{Z}_8[x]}{(x^2-4, 2x)}$& %
$\dfrac{\mathbb{Z}_4[x,y]}{(x^3,x^2-2,xy,y^2-2,y^3)}$& %
$\dfrac{\mathbb{Z}_4[x]}{(x^2)}$& %
$\,\,\,\dfrac{\mathbb{Z}_4[x,y]}{(x^2,y^2, xy-2)}\,\,\,$& %
$\dfrac{\mathbb{Z}_2[x,y]}{(x^2,y^2)}$&
$\mathbb{Z}_{27}$& %
$\dfrac{\mathbb{Z}_3[x]}{(x^3)}$& %
$\dfrac{\mathbb{Z}_9[x]}{(x^2-3,x^3)}$& %
$\dfrac{\mathbb{Z}_9[x]}{(x^2+3,x^3)}$& \\ %
&&&&&&&&&&\\
\hline
&&&&&&&&&&\\
$|R|$&16&16&16&16&16&27&27&27&27&\\
&&&&&&&&&&\\
$|R/\m|$&2&2&2&2&2&3&3&3&3&\\
&&&&&&&&&&\\
$char(R)$&8&4&4&4&2&27&3&9&9&\\
&&&&&&&&&&\\
$\Gamma(R)$&$G_3$&$G_3$&$G_4$&$G_4$&$G_4$&$G_1$&$G_1$&$G_1$&$G_1$&\\
&&&&&&&&&&\\
$|V(\Gamma(R))|$&7&7&7&7&7&8&8&8&8&\\
&&&&&&&&&&\\
$\gamma(\Gamma(R))$&0&0&0&0&0&0&0&0&0&\\
&&&&&&&&&&\\
\hline \multicolumn{11}{c}{\vspace{1mm}}\\
\multicolumn{11}{c}{\large Table 1. \,\,\,\, Local rings with planar zero-divisor graphs.}\\
\end{tabular}
\end{sideways}
\thispagestyle{empty} \pagebreak

\thispagestyle{empty} \tiny \centering
\begin{sideways}
\footnotesize $G_5$ is the graph mentioned in Remark~\ref{k34}.
\end{sideways}
\begin{sideways}
\begin{tabular}{|c|cccccccccc|}
\hline
&&&&&&&&&&\\
$R$& %
$\,\,\,\,\,\,\,\,\,\mathbb{Z}_{49}\,$& %
$\,\,\,\dfrac{\mathbb{Z}_7[x]}{(x^2)}\,\,\,$& %
$\dfrac{\mathbb{Z}_2[x,y]}{(x^3,xy,y^2)}$& %
$\dfrac{\mathbb{Z}_4[x]}{(x^3,2x)}$& %
$\dfrac{\mathbb{Z}_4[x,y]}{(x^3,x^2-2,xy,y^2)}$& %
$\dfrac{\mathbb{Z}_8[x]}{(x^2,2x)}$& %
$\dfrac{\mathbb{F}_8[x]}{(x^2)}$& %
$\dfrac{\mathbb{Z}_4[x]}{(x^3+x+1)}$& %
$\dfrac{\mathbb{Z}_4[x,y]}{(2x,2y, x^2,xy,y^2)}$& %
$\dfrac{\mathbb{Z}_2[x, y, z]}{(x, y, z)^2}$\\ %
&&&&&&&&&&\\
\hline
&&&&&&&&&&\\
$|R|$&49&49&16&16&16&16&64&64&16&16\\
&&&&&&&&&&\\
$|R/\m|$&7&7&2&2&2&2&8&8&2&2\\
&&&&&&&&&&\\
$char(R)$&49&7&2&4&4&8&2&4&4&2\\
&&&&&&&&&&\\
$\Gamma(R)$&$K_6$&$K_6$&$K_{1,1,1,4}$&$K_{1,1,1,4}$&$K_{1,1,1,4}$&$K_{1,1,1,4}$&$K_7$&$K_7$&$K_7$&$K_7$\\
&&&&&&&&&&\\
$|V(\Gamma(R))|$&6&6&7&7&7&7&7&7&7&7\\
&&&&&&&&&&\\
$\gamma(\Gamma(R))$&1&1&1&1&1&1&1&1&1&1\\
&&&&&&&&&&\\
\hline
&&&&&&&&&&\\
$R$& %
$\mathbb{Z}_{32}$& %
$\dfrac{\mathbb{Z}_2[x]}{(x^5)}$& %
$\dfrac{\mathbb{Z}_4[x]}{(x^3-2,x^5)}$& %
$\dfrac{\mathbb{Z}_4[x]}{(x^4-2, x^5)}$& %
$\dfrac{\mathbb{Z}_8[x]}{(x^2-2, x^5)}$& %
$\dfrac{\mathbb{Z}_8[x]}{(x^2-2x+2, x^5)}$& %
$\dfrac{\mathbb{Z}_8[x]}{(x^2+2x-2, x^5)}$&&&\\ %
&&&&&&&&&&\\
\hline
&&&&&&&&&&\\
$|R|$&32&32&32&32&32&32&32&&&\\
&&&&&&&&&&\\
$|R/\m|$&2&2&2&2&2&2&2&&&\\
&&&&&&&&&&\\
$char(R)$&32&2&4&4&8&8&8&&&\\
&&&&&&&&&&\\
$\Gamma(R)$&$G_5$&$G_5$&$G_5$&$G_5$&$G_5$&$G_5$&$G_5$&&&\\
&&&&&&&&&&\\
$|V(\Gamma(R))|$&15&15&15&15&15&15&15&&&\\
&&&&&&&&&&\\
$\gamma(\Gamma(R))$&1&1&1&1&1&1&1&&&\\
&&&&&&&&&&\\
\hline \multicolumn{11}{c}{\vspace{1mm}}\\
\multicolumn{11}{c}{\large Table 2. \,\,\,\, Local rings with toroidal zero-divisor graphs.}\\
\end{tabular}
\end{sideways}
\pagebreak

\thispagestyle{empty} \tiny \centering
\begin{sideways}
\footnotesize $p$ is prime and $p\neq 2,3$.
\end{sideways}
\begin{sideways}
\begin{tabular}{|c|cccccccccc|}
\hline
&&&&&&&&&&\\
$R$ &%
$\mathbb{Z}_2\times \mathbb{F}_{p^n}$& %
$\mathbb{Z}_3\times \mathbb{F}_{p^n}$& %
$\mathbb{Z}_2\times \mathbb{Z}_2$& %
$\mathbb{Z}_2\times \mathbb{Z}_3$& %
$\mathbb{Z}_2\times \mathbb{F}_4$& %
$\mathbb{Z}_3\times \mathbb{Z}_3$& %
$\mathbb{Z}_2\times \mathbb{Z}_4$& %
$\mathbb{Z}_2\times \dfrac{\mathbb{Z}_2[x]}{(x^2)}$& %
$\mathbb{Z}_3\times \mathbb{F}_4$& %
$\mathbb{Z}_3\times \mathbb{Z}_4$\\ %

&&&&&&&&&&\\
\hline
&&&&&&&&&&\\
$|R|$&$2p^n$&$3p^n$&4&6&8&9&8&8&12&12\\
&&&&&&&&&&\\
$|Spec(R)|$&2&2&2&2&2&2&2&2&2&2\\
&&&&&&&&&&\\
$char(R)$&$2p$&$3p$&2&6&2&3&4&2&6&12\\
&&&&&&&&&&\\
$|V(\Gamma(R))|$&$p^n$&$p^n+1$&2&3&4&4&5&5&5&7\\
&&&&&&&&&&\\
$\gamma(\Gamma(R))$&0&0&0&0&0&0&0&0&0&0\\
&&&&&&&&&&\\
\hline
&&&&&&&&&&\\
$R$& %
$\mathbb{Z}_3\times \dfrac{\mathbb{Z}_2[x]}{(x^2)}$& %
$\mathbb{Z}_2\times \mathbb{Z}_8$& %
$\mathbb{Z}_2\times \dfrac{\mathbb{Z}_2[x]}{(x^3)}$& %
$\mathbb{Z}_2\times \dfrac{\mathbb{Z}_4[x]}{(x^2-2, x^3)}$& %
$\mathbb{Z}_2\times \mathbb{Z}_9$& %
$\mathbb{Z}_2\times \dfrac{\mathbb{Z}_3[x]}{(x^2)}$& %
$\mathbb{Z}_3\times \mathbb{Z}_9$& %
$\mathbb{Z}_3\times \dfrac{\mathbb{Z}_3[x]}{(x^2)}$& %
$\mathbb{Z}_2\times \mathbb{Z}_2\times \mathbb{Z}_2$& %
$\mathbb{Z}_2\times \mathbb{Z}_2\times \mathbb{Z}_3$\\ %
&&&&&&&&&&\\
\hline
&&&&&&&&&&\\
$|R|$&12&16&16&16&18&18&27&27&8&12\\
&&&&&&&&&&\\
$|Spec(R)|$&2&2&2&2&2&2&2&2&3&3\\
&&&&&&&&&&\\
$char(R)$&6&8&2&4&18&6&9&3&2&6\\
&&&&&&&&&&\\
$|V(\Gamma(R))|$&7&11&11&11&11&11&14&14&6&9\\
&&&&&&&&&&\\
$\gamma(\Gamma(R))$&0&0&0&0&0&0&0&0&0&0\\
&&&&&&&&&&\\
\hline \multicolumn{11}{c}{\vspace{1mm}}\\
\multicolumn{11}{c}{\large Table 3. \,\,\,\,Non-local rings with planar zero-divisor graphs.}\\ %
\end{tabular}
\end{sideways}
\thispagestyle{empty} \pagebreak

\thispagestyle{empty} \tiny \centering
\begin{sideways}
\begin{tabular}{|c|cccccccccc|}
\hline
&&&&&&&&&&\\
$R$&
$\mathbb{F}_4\times \mathbb{F}_4$& %
$\mathbb{F}_4\times \mathbb{Z}_5$& %
$\mathbb{Z}_5\times \mathbb{Z}_5$& %
$\mathbb{Z}_4\times \mathbb{F}_4$& %
$\mathbb{F}_4\times \dfrac{\mathbb{Z}_2[x]}{(x^2)}$& %
$\mathbb{F}_4\times \mathbb{Z}_7$& %
$\mathbb{Z}_4\times \mathbb{Z}_4$& %
$\mathbb{Z}_4\times \dfrac{\mathbb{Z}_2[x]}{(x^2)}$& %
$\dfrac{\mathbb{Z}_2[x]}{(x^2)}\times \dfrac{\mathbb{Z}_2[x]}{(x^2)}$& %
$\mathbb{Z}_2\times \dfrac{\mathbb{Z}_2[x, y]}{(x^2, xy, y^2)}$\\ %

&&&&&&&&&&\\
\hline
&&&&&&&&&&\\
$|R|$&16&20&25&16&16&28&16&16&16&16\\
&&&&&&&&&&\\
$|Spec(R)|$&2&2&2&2&2&2&2&2&2&2\\
&&&&&&&&&&\\
$char(R)$&2&10&5&4&2&14&4&4&2&2\\
&&&&&&&&&&\\
$|V(\Gamma(R))|$&6&7&8&9&9&9&11&11&11&11\\
&&&&&&&&&&\\
$\gamma(\Gamma(R))$&1&1&1&1&1&1&1&1&1&1\\
&&&&&&&&&&\\
\hline
&&&&&&&&&&\\
$R$&
$\mathbb{Z}_2\times \dfrac{\mathbb{Z}_4[x]}{(2x, x^2)}$& %
$\mathbb{Z}_4\times \mathbb{Z}_5$& %
$\mathbb{Z}_5\times \dfrac{\mathbb{Z}_2[x]}{(x^2)}$& %
$\mathbb{Z}_3\times \mathbb{Z}_8$& %
$\mathbb{Z}_3\times \dfrac{\mathbb{Z}_2[x]}{(x^3)}$& %
$\mathbb{Z}_3\times \dfrac{\mathbb{Z}_4[x]}{(x^2-2, x^3)}$& %
$\mathbb{Z}_4\times \mathbb{Z}_7$& %
$\mathbb{Z}_7\times \dfrac{\mathbb{Z}_2[x]}{(x^2)}$& %
$\mathbb{Z}_2\times \dfrac{\mathbb{F}_4[x]}{(x^2)}$& %
$\mathbb{Z}_2\times \dfrac{\mathbb{Z}_4[x]}{(x^2+x+1)}$\\ %
&&&&&&&&&&\\
\hline
&&&&&&&&&&\\
$|R|$&16&20&20&24&24&24&28&28&32&32\\
&&&&&&&&&&\\
$|Spec(R)|$&2&2&2&2&2&2&2&2&2&2\\
&&&&&&&&&&\\
$char(R)$&4&20&10&24&6&12&28&14&2&4\\
&&&&&&&&&&\\
$|V(\Gamma(R))|$&11&11&11&15&15&15&15&15&19&19\\
&&&&&&&&&&\\
$\gamma(\Gamma(R))$&1&1&1&1&1&1&1&1&1&1\\
&&&&&&&&&&\\
\hline
&&&&&&&&&&\\
$R$& %
$\mathbb{Z}_2\times \mathbb{Z}_2\times \mathbb{F}_4$& %
$\mathbb{Z}_2\times \mathbb{Z}_2\times \mathbb{Z}_4$& %
$\mathbb{Z}_2\times \mathbb{Z}_2\times \dfrac{\mathbb{Z}_2[x]}{(x^2)}$& %
$\mathbb{Z}_2\times \mathbb{Z}_3\times \mathbb{Z}_3$& %
$\mathbb{Z}_2\times \mathbb{Z}_2\times \mathbb{Z}_5$& %
$\mathbb{Z}_2\times \mathbb{Z}_3\times \mathbb{F}_4$& %
$\mathbb{Z}_3\times \mathbb{Z}_3\times \mathbb{Z}_3$& %
$\mathbb{Z}_2\times \mathbb{Z}_2\times \mathbb{Z}_7$&
$\mathbb{Z}_2\times \mathbb{Z}_2\times\mathbb{Z}_2\times\mathbb{Z}_2$&\\ %
&&&&&&&&&&\\
\hline
&&&&&&&&&&\\
$|R|$&16&16&16&18&20&24&27&28&16&\\
&&&&&&&&&&\\
$|Spec(R)|$&3&3&3&3&3&3&3&3&4&\\
&&&&&&&&&&\\
$char(R)$&2&4&2&6&10&6&3&14&2&\\
&&&&&&&&&&\\
$|V(\Gamma(R))|$&12&13&13&13&15&17&18&21&14&\\
&&&&&&&&&&\\
$\gamma(\Gamma(R))$&1&1&1&1&1&1&1&1&1&\\
&&&&&&&&&&\\
\hline

\hline \multicolumn{11}{c}{\vspace{1mm}}\\
\multicolumn{11}{c}{\large Table 4. \,\,\,\,Non-local rings with
toroidal zero-divisor graphs.}\\
\end{tabular}
\end{sideways}

\thispagestyle{empty}

\end{document}